\documentclass[journal, letter]{IEEEtran}

\usepackage{amsmath}
\usepackage{amsfonts}
\usepackage{graphics}
\usepackage{pstricks}
\usepackage{array}
\usepackage{float}
\usepackage{epsf}
\usepackage{graphicx}
\usepackage{multirow}
\usepackage{cite}
\usepackage{acronym}
\usepackage{comment}
\usepackage{upgreek}
\usepackage{subcaption}
\usepackage{nicefrac} 
\usepackage[font=small]{caption}
    
\newcommand{\bfg}[1]{\boldsymbol{#1}}
\newcommand{\bfp}[1]{\boldsymbol{#1}'}
\newcommand{\bfpp}[1]{\boldsymbol{#1}''}
\newcommand{\bfppp}[1]{\boldsymbol{#1}'''}
\newcommand{\bfd}[1]{\dot{\boldsymbol{#1}}}
\newcommand{\bfdd}[1]{\ddot{\boldsymbol{#1}}}
\newcommand{\bfddd}[1]{\dddot{\boldsymbol{#1}}}
\newcommand{\bfb}[1]{\boldsymbol{\rm #1}}
\newcommand{\flux}{\bfg \varphi}
\newcommand{\ii}{\imath}
\newcommand{\jj}{\jmath}

\newcommand{\e}[1]{\boldsymbol{\rm e}_{\rm #1}}
\newcommand{\ep}[1]{\boldsymbol{\rm e}'_{\rm #1}}
\newcommand{\eh}[1]{\hat{\boldsymbol{\rm e}}_{\rm #1}}
\newcommand{\sw}{\bfg \omega}

\newcommand{\T}{\bfb T}
\newcommand{\B}{\bfb B}
\newcommand{\N}{\bfb N}
\newcommand{\Td}{\dot{\bfb T}}
\newcommand{\Bd}{\dot{\bfb B}}
\newcommand{\Nd}{\dot{\bfb N}}
\newcommand{\Tp}{\bfb T'}
\newcommand{\Bp}{\bfb B'}
\newcommand{\Np}{\bfb N'}
\newcommand{\wdq}{w_{\rm dq}}
\newcommand{\vd}{v_{\rm d}}
\newcommand{\vq}{v_{\rm q}}
\newcommand{\vo}{v_{\rm o}}
\newcommand{\vdp}{v'_{\rm d}}
\newcommand{\vqp}{v'_{\rm q}}
\newcommand{\vop}{v'_{\rm o}}
\newcommand{\vrp}{\hat{\bfg v}'}
\renewcommand{\th}[1]{\theta_{#1}}
\newcommand{\thp}[1]{\theta'_{#1}}
\newcommand{\V}[1]{V_{#1}}
\newcommand{\Vp}[1]{V'_{#1}}

\renewcommand{\sc}[2]{{\rm sc}(\th{#1}, \th{#2})}
\newcommand{\wo}{w_o}
\newcommand{\w}[1]{w_{#1}}

\acrodef{rocof}[RoCoF]{Rate of Change of Frequency}

\begin{document}

\title{Applications of the Frenet Frame to Electric Circuits}

\author{Federico~Milano,~\IEEEmembership{Fellow,~IEEE}, Georgios
  Tzounas,~\IEEEmembership{Member,~IEEE}, \\ Ioannis Dassios, and
  Taulant K\"{e}r\c{c}i, \IEEEmembership{Student~Member,~IEEE}%
  \thanks{F.~Milano, G.~Tzounas, I.~Dassios and T.~K\"{e}r\c{c}i are
    with School of Electrical and Electronic Engineering, University
    College Dublin, Dublin, D04V1W8, Ireland.  E-mails:
    \{federico.milano, ioannis.dassios, georgios.tzounas\}@ucd.ie}%
  \thanks{This work is supported by the European Commission by funding
    F.~Milano, G.~Tzounas and T.~K\"{e}r\c{c}i under project edgeFLEX,
    Grant No.~883710; and by Science Foundation Ireland by funding
    F.~Milano and I.~Dassios under project AMPSAS, Grant
    No.~SFI/15/IA/3074.}%
}

\maketitle

\begin{abstract}
  The paper discusses the relationships between electrical quantities,
  such as voltages, currents, and frequency, and geometrical ones,
  namely curvature and torsion.  The proposed approach is based on the
  Frenet frame utilized in differential geometry and provides a
  general framework for the definition of the time derivative of
  electrical quantities in stationary as well as transient conditions.
  As a byproduct, the proposed approach unifies and generalizes the
  time- and phasor-domain frameworks.  Other noteworthy results are a
  new interpretation of the link between frequency and the time
  derivatives of voltage and current; and a definition of the rate of
  change of frequency that includes the novel concept of ``torsional
  frequency.'' Several numerical examples based on balanced,
  unbalanced, harmonically-distorted and transient voltages illustrate
  the findings of the paper.
\end{abstract}


\IEEEpeerreviewmaketitle

\section{Notation}

In this paper,
scalars are indicated with normal
font, e.g.~$x$, whereas vectors are indicated in bold face,
e.g.~$\bfg x = (x_1, x_2, x_3)$.  All vectors have order 3, unless
otherwise indicated.

\vspace{3mm}

\subsubsection*{Scalars}

\begin{itemize}[\IEEEiedlabeljustifyl \IEEEsetlabelwidth{Z} \labelsep 0.9cm]
\item[$s$] length of a curve
\item[$t$] time
\item[$V$] voltage magnitude
\item[$w$] angular frequency 
\item[$\eta$] symmetric part of the geometric RoCoF
\item[$\theta$] voltage phase angle
\item[$\kappa$] curvature
\item[$\xi$] torsional frequency
\item[$\rho$] symmetric part of the geometric frequency
\item[$\tau$] torsion
\item[$\omega$] magnitude of vector $\bfg \omega$ 
\end{itemize}

\vspace{3mm}

\subsubsection*{Vectors}

\begin{itemize}[\IEEEiedlabeljustifyl \IEEEsetlabelwidth{Z} \labelsep 0.9cm]
\item[$\bfg 0$] null vector %
\item[$\B$] binormal vector of the Frenet frame %
\item[$\boldsymbol{\rm e}_i$] $i$-th vector of an orthonormal basis
\item[$\bfg \imath$] current vector
\item[$\bfg n$] normal vector before normalization
\item[$\N$] normal vector of the Frenet frame %
\item[$\bfg q$] electric charge vector
\item[$\T$] tangent vector of the Frenet frame %
\item[$\bfg v$] voltage vector 
\item[$\bfg \phi$] magnetic flux vector
\item[$\bfg \omega$] antisymmetric part of the geometric frequency
\end{itemize}

\vspace{3mm}

\subsubsection*{Derivatives}

\begin{itemize}[\IEEEiedlabeljustifyl \IEEEsetlabelwidth{Z} \labelsep 0.9cm]
\item[$x', \bfg x'$] derivative of a scalar/vector with respect to $t$
\item[$\dot{x}, \dot{\bfg x}$] derivative of a scalar/vector with respect to $s$
\item[$D^{\bfg x}_t$] time derivative operator applied to vector $\bfg x$
\end{itemize}

\section{Introduction}
\label{sec:intro}

\subsection{Motivation}

The study and simulation of circuit dynamics has traditionally been
approached using different frameworks.  Stationary AC circuits are
conveniently studied using quantities such as phasors and impedances;
circuits with harmonic contents are studied using Fourier analysis or
similar frequency-domain approaches; rotating machines and power
electronic devices are often studied using Park and/or Clarke
transforms; generic transients are studied using a time-domain
analysis \cite{FDF:2020}.  In this paper, we propose an approach based
on differential geometry, more specifically on the Frenet frame
\cite{Stoker}.  This approach leads to the definition of a framework
that admits, as special cases, the circuit analysis transformations
mentioned above.

\subsection{Literature Review}

Differential geometry finds applications in several fields of science
and engineering.  Some examples are the use of differential geometric
properties, such as that of curvature, in image segmentation and
three-dimensional object description \cite{1989yokoya}, as well as in
robotic control along geodesic paths \cite{2005selig}.  Another
relevant example are the utilities of the Frenet frame in the area of
autonomous vehicle driving \cite{2003lapierre,2010werling}.  Moreover,
there is a number of applications that are based on the theory of
geometric algebra, for example the use of quaternions in computer
graphics and visualization \cite{1995hanson, 2011vince} and in the
control of multi-agent networked systems \cite{2020talebi}.

The utilization of concepts of geometric algebra in circuit and power
system analysis is limited.  There is a group of works that elaborate
on the concept of instantaneous power \cite{4450060503, 1308315,
  4126799, 4512338, 5316097, IPT} that provide an interpretation of
the active and reactive power as the inner and cross (or wedge in the
polyphase case) products, respectively, of voltage and currents.  More
recently, some studies, including \cite{math9111295, 6161614,
  2016barry, 2018ishihara, 2015talebi}, have attempted to extend the
instantaneous power theory to a systematic study of electrical
quantities or circuits in the framework of geometric algebra.  In the
same vein, but using a novel perspective, \cite{freqgeom} makes an
additional step by proposing to interpret voltages and currents as the
\textit{time derivative of a multi-dimensional curve}.  This
interpretation allows the definition of the ``geometric frequency'' as
the result of an inner and an outer product.

In this work, we exploit differential geometry rather than geometric
algebra.  We are interested in the geometrical ``meaning'' of the time
derivative of electrical quantities such as voltage, current and
frequency.  With this aim, the formulas obtained in the paper are
deduced through the Frenet frame \cite{Stoker}.  The importance in
circuit analysis of the time derivatives of voltages and currents is
apparent as they are required in the constitutive equations of
capacitors and inductors.  The relevance of the \acf{rocof}, on the
other hand, is due to the increasing penetration, in the electric grid
all around the world, of renewable energy sources and the consequent
shift from synchronous to non-synchronous generation.  The \ac{rocof}
is, in turn, strictly related to the amount of available inertia in
the system \cite{8450880}.  The ability to estimate accurately the
\ac{rocof} is thus becoming an important aspect of the measurements
utilized by system operators.  As a matter of fact, several works
discuss the estimation of the \ac{rocof} from an instrumentation point
of view \cite{8253484, 8675542, 8630817, 9057685, 9162502}.

\subsection{Contributions}

We apply differential geometry to define a general framework for the
definition of electrical quantities and their time derivatives.  The
specific contributions of the paper are the following.
\begin{itemize}
\item The derivation of the expressions of the tangent, normal and
  binormal vectors of the Frenet frame in terms of the voltage (or
  current) of an electrical circuit.
\item A novel interpretation of the time derivative of any order of
  voltage and current in electrical circuits.
\item An expression of the \ac{rocof} which involves the definition of
  the novel concept of ``torsional frequency,'' which is also proposed
  and defined in the paper.
\item An example that shows that analytic signals commonly utilized in
  signal processing are a special case of the proposed framework in
  two dimensions.
\end{itemize}
The meaning and derivation of the vectors of the Frenet frame when
applied to electric quantities such as voltage, current and frequency
are duly discussed in the paper.

\subsection{Organization}

The remainder of the paper is organized as follows.  Section
\ref{sec:frenet} outlines the concepts of differential geometry that
are needed for the derivations of the theoretical results of this
work, which are given in Section \ref{sec:derivative}.  Section
\ref{sec:examples} illustrates the formulas of the time derivatives
through a series of examples.  The examples are aimed at showing that
the formulas derived in Section \ref{sec:derivative} admit as special
cases widely utilized frameworks such as DC circuits, phasors and Park
transform, as well as illustrate the formulas in unbalanced cases that
lead to the birth of time-variant curvature and torsion.  Section
\ref{sec:conc} draws conclusions and outlines future work.

\section{Frenet Frame of Space Curves}
\label{sec:frenet}

Let us consider a space curve
$\bfg x:[0,+\infty)\rightarrow\mathbb{R}^3$ with
$\bfg x = (x_1, x_2, x_3)$. Where $x_1=x_1(t)$, $x_2=x_2(t)$,
$x_3=x_3(t)$, is the set of parametric equations for the
curve. Equivalently:
\begin{equation}
    \bfg x = x_1 \, \e{1} + x_2 \, \e{2} + x_3 \, \e{3} \, ,
\end{equation}
where $(\e{1}, \e{2}, \e{3})$ is an orthonormal basis.  The length $s$
of the curve is defined as:
\begin{equation}
  s = \int_0^t \sqrt{\bfp x(r) \cdot \bfp x(r)} \, dr + s_0 \, ,
\end{equation}
from which one obtains the expression:
\begin{equation}
  \label{eq:s}
  s' = \frac{ds}{dt} = \sqrt{\bfp x \cdot \bfp x} = |\bfp x| \, ,
\end{equation}
where 
\begin{equation}
  \bfg x' = \frac{d}{dt} (x_1 \, \e{1}) +
  \frac{d}{dt} (x_2 \, \e{2}) + \frac{d}{dt} (x_3 \, \e{3}) \, ,
\end{equation}
and $\cdot$ represents the inner product of two vectors, which in
three dimensions, for $\bfg a = (a_1, a_2, a_3)$,
$\bfg b = (b_1, b_2, b_3)$, becomes:
\begin{equation}
  \label{eq:inner}
  \bfg a \cdot \bfg b = a_1b_1 + a_2b_2 + a_3b_3 \, .
\end{equation}
The length $s$ is an invariant of the curve.  It is relevant to
observe that, according to the chain rule, the derivative of $\bfg x$
with respect to $s$ can be written as:
\begin{equation}
  \bfd x = \frac{d \bfg x}{d s} =
  \frac{d\bfg x}{dt} \frac{dt}{ds} = \frac{\bfp x}{s'} =
  \frac{\bfp x}{|\bfp x|} \, .
\end{equation}
The vector $\bfd x$ has magnitude 1 and is tangent to the curve
$\bfg x$.

The Frenet frame is defined by the tangent vector $\T$, the normal
vector $\N$ and the binormal vector $\B$, as follows:
\begin{equation}
\label{eq:TNB}
  \begin{aligned}
    \T &= \bfd x \, , \\
    \N &= \frac{\bfdd x}{|\bfdd x|} \, , \\
    \B &= \T \times \N \, , \\
  \end{aligned}
\end{equation}
where $\times$ represents the cross product, which in three dimensions
can be written as the determinant of a matrix, as follows:
\begin{equation}
    \bfg a \times \bfg b = 
    \left | 
    \begin{matrix}
    \e{1} & \e{2} & \e{3} \\
    a_1 & a_2 & a_3 \\
    b_1 & b_2 & b_3 \\
    \end{matrix} 
    \right | \, .
\end{equation}
The vectors in \eqref{eq:TNB} are orthonormal,
i.e.~$\T = \N \times \B$ and $\N = \B \times \T$, and have relevant
properties, which can be expressed as follows \cite{Stoker}:
\begin{equation}
  \label{eq:frenet3}
  \begin{aligned}
    \Td &= \phantom{-}\kappa \N \, , \\
    \Nd &= -\kappa \T + \tau \B \, , \\
    \Bd &= - \tau \N \, ,
  \end{aligned}
\end{equation}
where $\kappa$ and $\tau$ are the \textit{curvature} and the
\textit{torsion}, respectively, which are given by:
\begin{equation}
  \label{eq:kappa}
  \kappa = |\bfdd x| = \frac{|\bfp x \times \bfpp x|}{|\bfp x|^3} \, , 
\end{equation}
and
\begin{equation}
  \label{eq:tau}
  \tau = \frac{\bfd x \cdot \bfdd x \times \bfddd x}{\kappa^2} =
  \frac{\bfp x \cdot \bfpp x \times \bfppp x}{|\bfp x \times \bfpp x|^2} \, .
\end{equation}
The quantities defined above, namely $\kappa$ and $\tau$, as well as
\eqref{eq:frenet3}, are utilized in the following section.

\section{Electrical Quantities in the Frenet Frame}
\label{sec:derivative}

This section presents the main theoretical results of the paper.  In
particular, the Frenet frame as well as of the curvature and torsion
of a space curve are expressed in terms of electrical quantities.
Then the expressions of the time derivatives of the vectors of
voltage, current as well as the frequency of these quantities are
derived based on the Frenet frame.  A general expression for
higher-order derivatives is also presented at the end of the section.

\subsection{Voltage and its Time Derivative in the Frenet Frame}
\label{sub:voltage}

The starting assumption of the discussion given in this section is
that the vector of the voltage, $\bfg v$, is the time derivative of a
space curve.  From a physical point of view, this means assuming that
the vector that describes the magnetic flux, say $\flux$, is formally
defined as:
\begin{equation}
  \label{eq:phi}
  \flux = - \bfg x \, .  
\end{equation}
Then, Faraday's law gives:
\begin{equation}
  \label{eq:v}
  \bfg v = -\bfp \flux = \bfp x \, .
\end{equation}
Then one can rewrite the expressions of the vectors $\T$, $\N$ and $\B$
of the Frenet frame in terms of the vector for the voltage and its
derivatives.

Let us observe first that the derivative of the length $s$, according
to \eqref{eq:s} and \eqref{eq:v}, becomes \cite{freqgeom}:
\begin{equation}
  \label{eq:sdot}
  s' = |\bfg v| = v \, ,
\end{equation}
and, then
\begin{equation}
  \label{eq:xd}
  \bfd x =
  -\bfd \flux =
  -\frac{\bfp \flux}{s'} =
  \frac{\bfg v}{v} \, ,
\end{equation}
and:
\begin{equation}
  \label{eq:xdd}
  \bfdd x =
  -\bfdd \flux =
  \frac{\bfp v}{v^2} -
  \frac{v' \, \bfg v}{v^3} \, ,
\end{equation}
and:
\begin{equation}
  \label{eq:xddd}
  \bfddd x =
  -\bfddd \flux =
  \frac{\bfpp v}{v^3} -
  3 \frac{v' \, \bfp v }{v^4} +
  3 \frac{(v')^2 \, \bfg v }{v^5} -
  \frac{v'' \, \bfg v}{v^4} \, ,
\end{equation}
where $v' = \frac{d}{dt}(v)$ and $v'' = \frac{d^2}{dt^2}(v)$.  It is
relevant to observe that, from the property
$\bfd x \cdot \bfdd x = 0$, as these vectors are orthogonal by
construction, and from \eqref{eq:xd} and \eqref{eq:xdd}, one obtains
\cite{freqgeom}:
\begin{equation}
  \label{eq:rho}
  \rho = \frac{v'}{v} \, .
\end{equation}
As it is well known, in time-frequency analysis and signal processing,
the quantity $\rho$ is defined as the \textit{instantaneous bandwidth}
\cite{Cohen:1995}.  In this work, however, we rather use the
interpretation of $\rho$ given in \cite{freqgeom}, namely, the
symmetric part of the \textit{geometric frequency}.\footnote{In this
  work, the terms symmetric and antisymmetric do not refer to the
  properties of a matrix but rather to the effect of operators.}  It
is also relevant to note that, from a geometrical point of view,
$\rho \, v = v'$ can be viewed as the ``radial'' component of the
velocity $\bfg v$.  In this vein, $\rho$ can be defined as
\textit{radial frequency}.

On the other hand, from \eqref{eq:kappa}, \eqref{eq:v} and
\eqref{eq:xd}-\eqref{eq:rho}, one has:
\begin{equation}
  \label{eq:kappa2}
  \boxed{\kappa = \frac{|\bfg v \times \bfp v|}{v^3} =
  \frac{|\sw|}{v}  = \frac{\omega}{v}} \, ,
\end{equation}
where the vector $\bfg \omega$ is defined as the antisymmetric
component of the geometric frequency, as follows
\cite{freqgeom}:\footnote{ Note that in \cite{freqgeom} the wedge
  product is used rather than the cross product, as the definition of
  the geometric frequency is given for a voltage vector of arbitrary
  dimension $n$.  However, for simplicity but without lack of
  generality, this paper focuses on the cases for which $n \le 3$.}
\begin{equation}
  \label{eq:rhoOmega}
  \bfg \omega = \frac{\bfg v \times \bfp v}{v^2}  \, .
\end{equation}
From a geometrical point of view, in 3 dimensions,
$\omega \, v = |\bfg \omega| \, v$ can be interpreted as the azimuthal
component of the velocity $\bfg v$.  Then, $\omega$ can be defined as
\textit{azimuthal frequency}.

Then, using the definition of $\bfg \omega$ above, the torsion given
in \eqref{eq:tau} can be rewritten as:
\begin{equation}
  \boxed{\tau = \frac{\bfg v \cdot \bfp v \times \bfpp v}{ \omega^2 \, v^4}} \, .
\end{equation}
The vectors of the Frenet frame can be written as:
\begin{equation}
  \label{eq:frenet2}
  \begin{aligned}
    \T = \frac{\bfg v}{v} \, , \qquad
    \N = \frac{\bfg n}{n} \, , \qquad
    \B = \frac{\sw}{\omega} \, ,
  \end{aligned}
\end{equation}
where $\bfg n$ is the normal vector before normalization, as follows:
\begin{equation}
  \label{eq:n}
  \begin{aligned}
    \bfg n &= \bfp v - \rho \bfg v \, , \\
    n &= |\bfg n| = \sqrt{|\bfp v|^2 - (\rho v)^2} \, .
  \end{aligned}
\end{equation}

Note that, from the following property of the scalar triple product:
\begin{equation}
  \label{eq:triple1}
  \bfg a \cdot \bfg b \times \bfg c = \bfg c \cdot \bfg a \times \bfg b \, , 
\end{equation}
the expression of the torsion can be rewritten as follows:
\begin{equation}
  \label{eq:tau2}
  \tau = \frac{\bfpp v \cdot \bfg v \times \bfp v}{\omega^2 v^4} 
  = \frac{\bfpp v \cdot \sw}{\kappa^2} \, ,
\end{equation}
which indicates that the torsion is null, apart from the obvious cases
$\bfpp v = \bfg 0$ and $\sw = \bfg 0$, if $\bfpp v$ is perpendicular
to $\sw$.  This happens if the voltage vector is unbalanced, as
illustrated in Section \ref{sec:examples}.

We are now ready to present one of the main results of this paper.  
Recalling that the Frenet vectors are orthonormal and, in particular,
$\N = \B \times \T$, one has:
\begin{equation}
  \frac{\bfg n}{n} = \frac{\sw}{\omega} \times \frac{\bfg v}{v} \, .
\end{equation}
Noting that $n$ is equal to the azimuthal speed, i.e.~$n=\omega v$
(see the proof in the Appendix), the expression above can be
simplified as:
\begin{equation}
  \bfg n = \sw \times \bfg v \, ,
\end{equation}
and, from \eqref{eq:n}:
\begin{equation}
  \label{eq:vdot}
  \boxed{\bfp v = \rho \, \bfg v + \sw \times \bfg v} \, .
\end{equation}
The previous expression is relevant as it allows the definition of a
time derivative operator.  Equation \eqref{eq:vdot} can be, in fact,
written as:
\begin{equation}
  \label{eq:vdot2}
  \left (\frac{d}{dt} - [\rho + \sw \times] \right )\bfg v=\bfg 0 \, ,
\end{equation}
which leads to define the operator:
\begin{equation}
  \label{eq:ddt}
  D_t^{\bfg v}:=\frac{d}{dt} - [\rho + \sw \times] \, .
\end{equation}
In \eqref{eq:ddt}, the upper symbol $\bfg x$ in the notation
$D^{\bfg x}_t$ indicates the vector upon which the operator is
applied, whereas the subindex $t$ indicates the independent variable
with respect to which the operator is calculated, i.e.~time.

Equation \eqref{eq:vdot} and the operator \eqref{eq:ddt} can be
interpreted as follows: the time derivative of a vector can be split
into two components, one symmetric ($\rho$) and one antisymmetric
($\sw \times$).  This operator has been obtained without any
assumption on the time dependence of $\bfg v$ nor on its dimension.
For dimensions greater than 3, in fact, it is sufficient to substitute
the cross product with the wedge product \cite{freqgeom}.

While it has been obtained considering the vector of the voltage, the
time derivative operator defined in \eqref{eq:ddt} can be applied to
any vector that represents the first time derivative of a space curve.
In particular, if one considers the vector of the current
$\bfg \ii = \bfp q$, being $\bfg q$ the vector of the electric charge
in a given point of a circuit \cite{freqgeom}, then:
\begin{equation}
  \label{eq:idot}
  D_t^{\bfg \ii} \, \bfg \ii = \bfg 0 \, .
\end{equation}
For simplicity, in the remainder of this work, we focus exclusively on
the voltage.  However, all examples provided in Section
\ref{sec:examples} can be equivalently applied to currents.

\subsection{Derivative of the Vector Frequency in the Frenet Frame}
\label{sub:rocof}

To obtain an expression of $\sw'$, we start with \eqref{eq:frenet3}
and rewrite the vectors in terms of the derivative with respect to
time, recalling that, for the chain rule:
\begin{equation}
  \Tp = v \Td \, , \qquad \Np = v \Nd \, , \qquad \Bp = v \Bd \, .
\end{equation}
Then, one obtains:
\begin{equation}
  \label{eq:frenet4}
  \begin{aligned}
    \Tp &= \phantom{-}\omega \N \, , \\
    \Np &= -\omega \T + \xi \B \, , \\
    \Bp &= - \xi \N \, ,
  \end{aligned}
\end{equation}
where $\omega = v \kappa$ from \eqref{eq:kappa2} and: 
\begin{equation}
  \label{eq:xi}
  \xi = v \tau \, ,
\end{equation}
where $\xi$ can be defined as the \textit{torsional frequency} and has
the unit of s$^{-1}$, as $\omega$ and $\rho$.  From \eqref{eq:frenet2}
and letting $\omega' =\frac{d}{dt}\omega$, the third equation of
\eqref{eq:frenet4} can be rewritten as follows:
\begin{equation}
  \frac{d}{dt} \frac{\sw}{\omega} =
  \frac{\sw'}{\omega} - \frac{\sw \omega'}{\omega^2} =
  -\xi \frac{\bfg n}{n} \, ,
\end{equation}
or, equivalently,
\begin{equation}
  \label{eq:rocof}
  \sw' = \omega'\frac{\sw}{\omega} -
  \omega \xi \frac{\bfg n}{n} \, ,
\end{equation}
or, equivalently,
\begin{equation}
\label{eq:rocof2}
\boxed{\sw' = \eta \, \sw + \tau \, \bfg v \times \sw} \, ,
\end{equation}
where $\eta = \omega'/\omega$ and $[\tau \bfg v \times]$ are the
symmetric and antisymmetric parts, respectively, of the time
derivative of $\sw$.

Equation \eqref{eq:rocof2} expresses the generalization of the
\ac{rocof}, commonly utilized in power system studies as a metric of
the severity of a transient following a contingency or a power
imbalance.  For a balanced system, $\tau$ is null and, thus,
\eqref{eq:rocof} leads to:
\begin{equation}
  |\sw'| = \omega' \, ,
\end{equation}
which represents the conventional definition of \ac{rocof}.  However,
\eqref{eq:rocof} shows that, when the torsion is not null, e.g., in
unbalanced cases, $\sw'$ is a vector with richer information than the
usual understanding of the \ac{rocof}.  Section \ref{sec:examples}
illustrates \eqref{eq:rocof} through numerical examples.

Finally, the time derivative of the symmetric part of the geometric
frequency defined in \cite{freqgeom} is:
\begin{equation}
  \label{eq:rhodot}
  \rho' = \frac{\bfg v \cdot \bfpp v}{v^2} +
  \omega^2 - \rho^2 = \frac{v''}{v} - \rho^2 \, , 
\end{equation}
which is obtained from the first equation of \eqref{eq:rhoOmega} and
\eqref{eq:rho}.

\subsection{Higher-Order Time Derivatives}

The vectors of the Frenet frame constitute a basis for 3-dimensional
systems, hence any vector, including any time derivatives of the
voltage, current, and frequency can be written in terms of $\N$, $\T$
and $\B$.  This can be deduced from \eqref{eq:frenet4}.  For example,
the second time derivative of $\T$ becomes:
\begin{equation}
  \T'' = \omega' \, \N + \omega \, \N' = \omega' \, \N -
  \omega^2 \, \T + \omega \xi \, \B \, ,
\end{equation}
or, equivalently, in terms of $(\bfg v, \bfg n, \sw)$:
\begin{equation}
  \bfpp v = (\rho' + \rho^2 - \omega^2) \, \bfg v -
  (2\rho - \eta) \, \bfg n - v \xi \, \sw \, .
\end{equation}
While the complexity of the expressions of the coefficients of
higher-order derivatives increases, the following general expression
for the $r$-th derivative holds:
\begin{equation}
  \label{eq:vdotr}
  \boxed{\bfg v^{(r)} = a_r \, \bfg v + b_r \, \bfg n + c_r \, \sw} \, .
\end{equation}
For example, for the first two derivatives, one has:
\begin{align*}
  a_1 &= \rho \, ,& b_1 &= 1 \, ,& c_1 &= 0 \, , \\
  a_2 &= \rho' + \rho^2 - \omega^2 \, ,& b_2 &= 2\rho - \eta \, ,& c_2 &= v\xi \, .
\end{align*}
The following remarks are relevant.
\begin{itemize}
\item Expressions with the same structure as \eqref{eq:vdotr} can be
  written also for $\bfg n^{(r)}$ and $\sw^{(r)}$.
\item Since $\bfg n = \sw \times \bfg v$, equation \eqref{eq:vdotr}
  is, in turn, a function exclusively of $\bfg v$ and $\sw$.
\end{itemize}

\section{Examples}
\label{sec:examples}

This section illustrates the theoretical results above through special
cases that are relevant in circuit and power system analysis.  The
first two examples are aimed at illustrating the properties of
\eqref{eq:ddt} in the special cases of DC circuits and stationary AC
circuits.  The subsequent examples show the effect of imbalances and
harmonics on the various components of the frequency of three-phase
voltages, as well as on the time derivative of the frequency itself.
In all examples below, we assume that voltages are curves in three
dimensions, where the basis, unless otherwise stated, is given by the
following orthonormal vectors:
\begin{equation}
  \label{eq:basis0}
  \begin{aligned}
    \e{1} &= (1, 0, 0) \, , \\ 
    \e{2} &= (0, 1, 0) \, , \\
    \e{3} &= (0, 0, 1) \, . 
  \end{aligned}
\end{equation}
As already said in the theoretical sections, systems with dimensions
higher than 3 can be also considered by using the wedge product rather
than the cross product and the generalization of the Frenet-Serret
formulas.  However, the study of voltages with dimensions higher than
three is beyond the scope of this paper.

In Sections~\ref{sub:3ph} to \ref{sub:variable}, we show state-space
three-dimensional plots of the voltage $\bfg v$ as this quantity has a
clear physical meaning and it is widely used in practice.  However, it
is important to keep in mind that, in the proposed framework, the
voltage is, in effect, the time derivative of a position (flux) vector
$\bfg x = -\bfg \varphi$.  The curvature and torsion, thus, are those
of such a flux vector, not of the voltage.

In the following, we do not discuss explicitly the behavior of the
vectors $\bfb T$, $\bfb N$ and $\bfb B$.  However, it is relevant to
observe that the voltage vector $\bfg v$ and the antisymmetric part of
the geometric frequency $\bfg \omega$, are the tangent and binormal
vectors, namely $\bfb T$ and $\bfb B$, of the Frenet frame before
normalization.  Moreover, $\bfg n$, which is the normal vector
$\bfb N$ of the Frenet frame before normalization, has an important
role in the definition of the time derivative of the voltage,
$\bfg v'$, as it is the antisymmetric part of such a derivative.
Thus, the discussions on the behavior of $\bfg v$, $\bfg v'$ and
$\bfg \omega$ in the examples below are, indirectly, also discussions
of the behavior of the vectors of the Frenet frame.

\subsection{DC Voltage}
\label{sub:dc}

As a first example, we show the effect of \eqref{eq:ddt} on DC
voltages.  From the geometrical point of view, a DC voltage is
equivalent to the time derivative of a curve with one dimension,
i.e. a straight line.  Using the vector notation, in DC, the voltage
is a curve that has only one component along one direction, say
$\e{1}$ of the basis.  Then:
\begin{equation}
  \label{eq:vdc}
  \bfg v = v_{\rm dc} \, \e{1} + 0 \, \e{2} + 0 \, \e{3} \, .
\end{equation}
It is immediate to show that, for a curve in one dimension,
$\kappa = \tau = 0$, $\sw = \bfg 0$ and, the operator \eqref{eq:ddt}
is reformulated as:
\begin{equation}
  \label{eq:rhodc}
  D^{\bfg v}_t = \frac{d}{dt} - \rho \, .
\end{equation}
Thinking in terms of curves, this result comes with no surprise, as a
straight line cannot rotate or twist.  Only the radial component of
the velocity, thus, can be nonnull.  It is also relevant to note that,
according to \eqref{eq:ddt}, expression \eqref{eq:rhodc} states that
the antisymmetric component of the time derivative of DC quantities is
always null.

\subsection{Stationary Single-Phase AC Voltage}
\label{sub:1ph}

Let us consider a stationary single-phase voltage with constant
angular frequency $\wo$ and magnitude $V$.  Then the voltage vector
can be written as:
\begin{equation}
\label{eq:single}
    \bfg v = V \cos(\wo t + \alpha) \, \e{1} + 
    V \sin(\wo t + \alpha) \, \e{2} + 0 \, \e{3} \, ,
\end{equation}
where $\alpha$ is a constant phase shift. 

Note that the representation of the voltage in \eqref{eq:single} is
that of an \textit{analytic signal}, that is, the first component is
the signal itself, namely $u(t) = V \cos(\wo t + \alpha)$, whereas the
second component is the Hilbert transform of the signal, i.e.,
$\mathcal{H}[u(t)] = V \sin(\wo t + \alpha)$ \cite{Hahn:1996,
  4026700}.  It is important to note that analytic signals are defined
as complex quantities rather than vectors but, nonetheless, they are
intrinsically two dimensional quantities.  The utilization of the
Frenet framework and the interpretation of the signal as a
\textit{curve} is more general and, as shown in this example, admits
analytic signals as a special case.

It is also relevant to note that, a single phase AC voltage
represents, from the geometric point of view, a plane curve.  This is
also consistent with the fact that an AC signal requires two
independent quantities to be defined, e.g., magnitude and phase angle
and, hence, the coordinate basis requires two dimensions to be
complete.
Applying the definitions given in the previous section, one can easily
find that:
\begin{align}
  \label{eq:rhoac}
  \begin{aligned}
  \rho &= 0 \, ,& \qquad
  \sw &= \wo \, \e{3} \, , \qquad \xi = 0 \, .
  \end{aligned}
\end{align}
These results were expected for a plane curve, which can rotate
($\kappa = \wo/V \ne 0$) but cannot twist ($\tau = \xi/V = 0$).  On
the other hand, $\rho = 0$ is a consequence of the fact that
$V = \rm const$.  Moreover, from \eqref{eq:ddt}, one obtains:
\begin{equation}
  \label{eq:ddtac}
  D^{\bfg v}_t = \frac{d}{dt} - \wo \, [\e{3} \times ] \, .
\end{equation}
It is worth to further elaborate on \eqref{eq:ddtac}.  We note first
that, since $\rho$, $\kappa$ and $\tau$ are geometric invariants, same
results can be obtained using any other orthonormal basis.  In
particular, if one chooses a basis that rotates at constant angular
speed $\wo$:
\begin{equation}
  \label{eq:basisac}
  \begin{aligned}
  \eh{1} &= (\cos(\wo t), 0, 0) \, , \\ 
  \eh{2} &= (0, \sin(\wo t), 0) \, , \\
  \eh{3} &= (0, 0, 1) \, , 
  \end{aligned}
\end{equation}
then the voltage vector becomes:
\begin{equation*}
  \bfg v =
  V \cos(\alpha) \, \eh{1} + V \sin(\alpha) \, \eh{2} =
  v_1 \, \eh{1} + v_2 \, \eh{2} \, ,
\end{equation*}
with $v = |\bfg v| = V$ and: 
\begin{equation}
  \label{eq:vdotac}
  \bfp v = \wo \, \eh{3} \times \bfg v = 
  - \wo v_2 \, \eh{1} + \wo v_1 \, \eh{2} \, ,
\end{equation}
where we have used the identities:
\begin{equation*}
  -\eh{1} = \eh{3} \times \eh{2} \, ,
  \qquad \eh{2} = \eh{3} \times \eh{1} \, .
\end{equation*}
Equation \eqref{eq:vdotac} has a striking formal similarity with the
well-known expression in phasor-domain where the derivative is given
by $\jj \, \wo$, where $\jj$ is the imaginary unit.  As a matter of
fact, in $\mathbb{R}^2$, the cross (wedge) vector is isomorphic to
complex numbers.  Thus, one can define the following correspondences:
\begin{equation}
  \begin{aligned}
    \eh{1} &\Rightarrow 1 \, , \\
    \eh{2} &\Rightarrow \jj \, , \\
    \wo [\eh{3} \times] &\Rightarrow \jj \, \wo \, ,
  \end{aligned}
\end{equation}
which leads to rewrite the vector $\bfg v$ as the well-known phasors,
namely $\bar v = v_1 + \jj \, v_2$.
We note also that, in two dimensions, $|\sw| = \wo$ hence the
azimuthal frequency coincides with the well-known angular frequency of
circuit analysis and with the \textit{instantaneous frequency} as
commonly defined in time-frequency analysis and signal processing
based on analytic signals \cite{Cohen:1995, Hahn:1996}.

Noteworthy, this observation can be generalized for any signal $u(t)$
for which the Hilbert transforms $\mathcal{H}[u(t)]$ exists.  With
this assumption, the analytic signal associated with $u(t)$ is the
following complex quantity:
\begin{equation}
  \bar{u}(t) = u(t) + j \mathcal{H}[u(t)] = u(t) + j \hat{u}(t) \, ,
\end{equation}
the instantaneous frequency of which is defined as \cite{Cohen:1995}:
\begin{equation}
  \label{eq:phidot}
  \phi'(t) = \frac{\hat{u}'(t)u(t) - u'(t)\hat u(t)}{u^2(t) + \hat u^2(t)} \, ,
\end{equation}
where 
\begin{equation}
    \phi (t) = {\rm arctan} \left ( \frac{\hat{u}(t)}{u(t)} \right )
\end{equation}
is the phase angle of $\bar{u}(t)$.  Using the approach proposed in
this work, on the other hand, we define the vector:
\begin{equation}
  \bfg u (t) = u(t) \, \e{1} + \hat{u}(t) \, \e{2}  + 0 \, \e{3}\, ,
\end{equation}
which leads to:
\begin{align}
  \label{eq:u}
  \begin{aligned}
  \rho &= \frac{u'u + \hat{u}'\hat{u}}{u^2 + \hat{u}^2} \, ,  &\quad
  \sw &= \frac{\hat{u}'u - u'\hat u}{u^2 + \hat u^2} \, \e{3} \, , &\quad \xi &= 0 \, ,
  \end{aligned}
\end{align}
where the time dependency has been omitted for simplicity.  Both
formulations, thus, leads to the same expression for the instantaneous
frequency, i.e., $\phi' = |\sw|$.  However, the proposed vector-based
approach is more general as it allows defining the quantities $\rho$
and $\xi$ and is not limited to two dimensions.

\subsection{Three-Phase AC Voltages}
\label{sub:3ph}

In this section, we consider a three-phase AC system.  A possible
approach is to utilize the same coordinates that we have considered
for the single-phase AC voltage of the previous example.  This
approach is the one commonly utilized in circuit analysis.  In this
section, however, we show that choosing as coordinates the voltages of
each phase leads to interesting results.  Hence, we define the voltage
vector as:
\begin{equation}
    \bfg v = v_a \, \e{1} + v_b \, \e{2} + v_c \e{3} \, .
    \label{eq:3ph:vector}
\end{equation}
Let us first consider the case of voltages represented by a single
harmonic:
\begin{equation}
  \begin{aligned}
    \label{eq:3ph:general}
    v_a &= \V{a} \, \sin (\th{a}) \, , \\
    v_b &= \V{b} \, \sin (\th{b}) \, , \\
    v_c &= \V{c} \, \sin (\th{c}) \, .
  \end{aligned}
\end{equation}
The expressions of $\rho$, $\sw$ and $\xi$ for \eqref{eq:3ph:general}
are:
\begin{equation}
  \label{eq:3ph:rho}
  \rho = \frac{\sum_{i} \V{i}^2 \thp{i}\sin(2\th{i}) +
    \V{i}\Vp{i}(1-\cos(2\th{i}))}{v^2} \, ,
\end{equation}
\begin{equation}
  \label{eq:3ph:omega}
  \sw = \frac{\sum_{ijk} (r_{jk} + u_{jk}) \, \bfb e_i}{v^2} \, ,
\end{equation}
\begin{equation}
  \label{eq:3ph:xi}
  \xi = \frac{v\,\sum_i (p_i \sin(\th{i}) + q_i \cos(\th{i})) \,
    \omega_i}{\sum_{jk} (r_{jk} + u_{jk})^2} \, ,
\end{equation}
where $i \in \{a, b, c\}$, $ijk \in \{abc, bca, cab\}$,
$jk \in \{ bc, ca, ab\}$ and:
\begin{align*}
  v &= \sqrt{\textstyle \sum_{i}\V{i}^2 (1 - \cos(2\th{i}))} \, , \\
  r_{jk} &= (\V{j}\Vp{k} - \V{k}\Vp{j}) \sin (\th{j}) \sin (\th{k}) \, , \\
  u_{jk} &= \V{j}\V{k} (\thp{k} \sin (\th{j}) \cos (\th{k}) -
           \thp{j} \sin (\th{k}) \cos (\th{j})) \, , \\
  p_i &= \V{i}\V{i}'' + (\Vp{i})^2 - \V{i}(\thp{i})^2 \, , \\
  q_i &= \Vp{i}\thp{i} - \V{i}\th{i}'' \, .
\end{align*}

The following remarks are relevant.
\paragraph*{Remark 1} In general, $\rho$, $\sw$ and $\xi$ depend on
the time derivatives of both the magnitudes and the phase angles of
the three-phase voltages.  This result is counter-intuitive.  In the
common understanding, in fact, the angular frequency is defined as the
time derivative of the sole phase angle of the voltage (see
\cite{IEEE118}).  The conventional expression of the angular frequency
is obtained if $\Vp{i} = 0$.  This result reiterates that the common
definition of ``frequency'' is, in effect, a special case of the
framework proposed in this paper.

\paragraph*{Remark 2} It is possible to have $\rho = 0$ and
$\sw \ne \bfg 0$ (for example, the obvious case of balanced and
stationary three-phase voltages, which is illustrated below) but also
$\rho \ne 0$ and $\sw = \bfg 0$ (for example a balanced voltage with
$\V{i}' \ne 0$).

\paragraph*{Remark 3} The torsional frequency $\xi$ is non-null if and
only if $\sw\ne \bfg 0$.  This is consequence of the expression of the
torsion given in \eqref{eq:tau2}.  The condition $\xi \ne 0$ holds for
$\V{i}'' \ne 0$ and/or $\th{i}'' \ne 0$.

\paragraph*{Remark 4} A question might arise on why the base
considered in Section \ref{sub:1ph} for the single-phase AC system is
not also utilized for the three-phase system.  One can, of course,
consider each phase of the three-phase system separately, effectively
considering each phase as a single-phase system as in the previous
example.  This is the common approach in three-phase circuit analysis
based on phasors.  However, in this and in the following examples, the
three-phase voltages are assumed to form a three-dimensional vector.
It is relevant to note that, since $\rho$, $\kappa$ and $\tau$ are
geometric invariants, it does not matter which coordinates one chooses
as long as these coordinates form a complete basis for the system.

\begin{figure*}[ht!]
  \centering
  \resizebox{\linewidth}{!}{\includegraphics{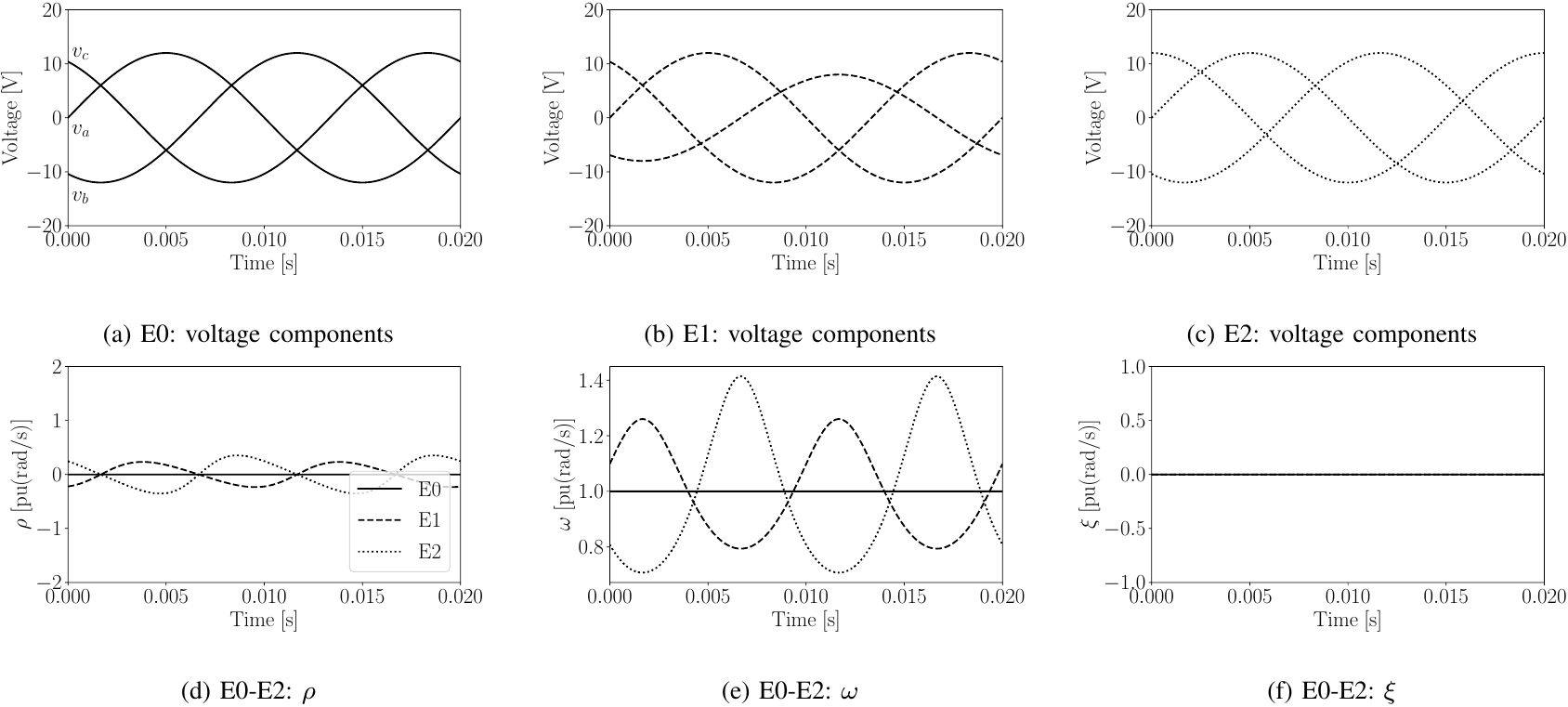}}
  \caption{E0-E2: Three-phase AC voltage components, and geometric
    invariants $\rho$, $\omega$ and $\xi$.}
  \label{fig:E0E1E2}
\end{figure*}
\begin{figure}[ht!]
  \centering
  \resizebox{0.7\linewidth}{!}{\includegraphics{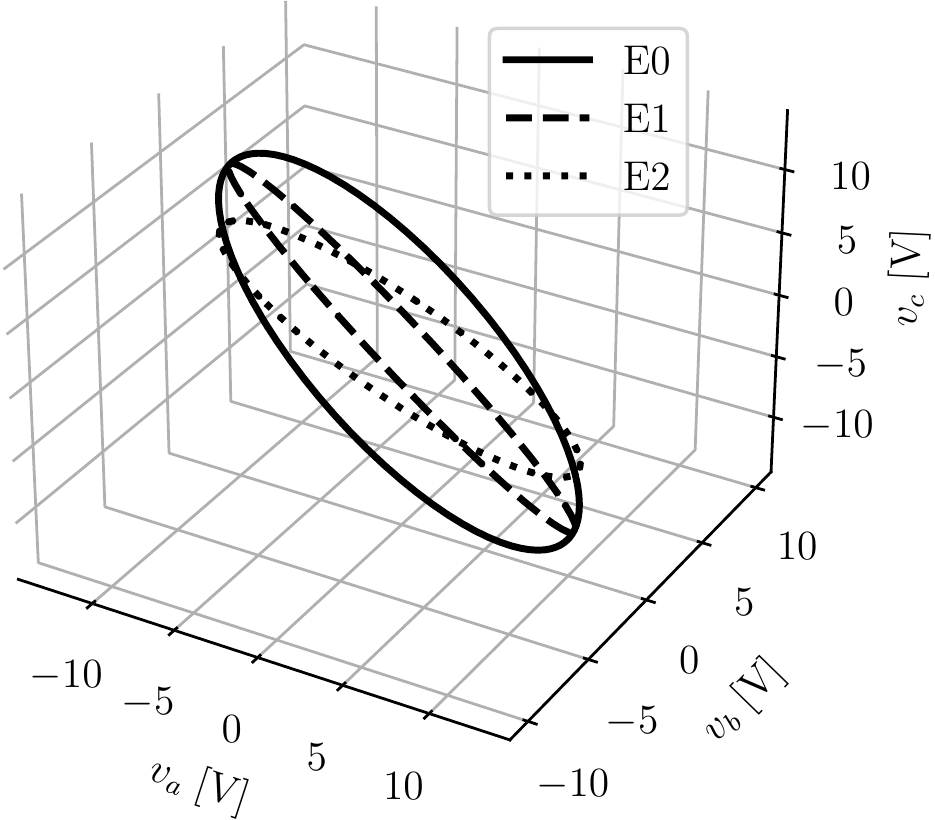}}
  \caption{Three-phase voltage in the space $(v_a, v_b, v_c)$, E0-E2.}
  \label{fig:vavbvc:E0E1E2}
\end{figure}

\paragraph*{Stationary Voltages}

Assuming that the voltage is stationary with constant angular
frequency $\wo$, the three components of the vector are:
\begin{equation}
  \begin{aligned}
  \label{eq:3ph:voltage}
    v_a &= V_a \sin (\wo t + \th{ao})
   \, , \\ 
   v_b &= V_b \sin (\wo t + \th{bo})
   \, , \\ 
   v_c &= V_c \sin (\wo t + \th{co})
   \, ,
  \end{aligned}
\end{equation} 
where $\V{i}$ is the voltage magnitude of phase $i$; and $\th{io}$
denotes the voltage angle of phase $i$ at $t=0$~s.
Let $\wo=100\pi$~rad/s, $\theta_{ao}=0$~rad.  Next, we consider
relevant special cases of \eqref{eq:3ph:general}.
%
\paragraph*{Positive and negative sequence voltages}
From \eqref{eq:3ph:rho}-\eqref{eq:3ph:xi}, one has
for the stationary positive sequence with $\V{a}=
\V{b}=\V{c} = \rm const.$ and $\th{bo} = -\th{co}= - {2\pi}/{3}$ rad:
\begin{equation*}
  \rho = \xi = 0 \, , \quad \sw = \frac{\wo}{\sqrt{3}}
  (\e{1} + \e{2} + \e{3}) \, .
\end{equation*}
    
Analogously, the stationary negative sequence, namely
$\V{a}= \V{b}=\V{c} = \rm const.$ and $\th{bo} = -\th{co}= {2\pi}/{3}$
rad, leads to:
\begin{equation*}
  \rho = \xi = 0 \, , \quad \sw = -\frac{\wo}{\sqrt{3}}
  (\e{1} + \e{2} + \e{3}) \, .
\end{equation*}
Hence, for the stationary positive and negative sequences, the time
derivative operator \eqref{eq:ddt} becomes:
\begin{equation}
  D^{\bfg v}_t =  \frac{d}{dt} \mp \frac{\wo}{\sqrt 3}
  [(\e{1} + \e{2} + \e{3}) \times] \, ,
\end{equation}
and, as in the example of Section V.B, $\omega = |\sw| = \wo$, namely
the azimuthal and angular frequencies coincide.

\paragraph*{Unbalanced voltage}
Unbalanced voltages are characterized by $\V{a} \ne \V{b} \ne \V{c}$
and/or $\th{bo} \ne - \th{co} \ne -2\pi/3$ rad.  Then, assuming that
the voltage magnitudes and phase shifts are constants, from
\eqref{eq:3ph:rho}-\eqref{eq:3ph:xi}, one has:
\begin{equation*}
  \begin{aligned}
    \rho &= \frac{\wo \sum_i \V{i}^2 \sin(2\th{i})}{v^2} \, ,
    \quad &\xi = 0 \, , \\ \sw &=
    \frac{\wo \sum_{ijk} \V{j}\V{k} \sin(\th{j} - \th{k}) \e{i}}{v^2} \, .
  \end{aligned}
\end{equation*}
In this case, thus, the components of the time derivative in
\eqref{eq:vdot} depend on the parameters of the voltage.
    
\paragraph*{Zero-sequence voltage}
This example allows us discussing the issue of the choice of the basis
for the voltage vector.
A zero-sequence voltage is composed of three equal AC voltages.  In
this case, thus, we cannot utilize the same coordinates we have used
in the previous two examples, namely $(v_a, v_b, v_c)$, as these are
linearly dependent and do not form a basis.  One has to proceed as
discussed in Section \ref{sub:1ph} for a single-phase AC voltage, for
which $\sw = \wo \, \e{3}$.  Thus, when considering three-phase
voltages, one should thus first remove the zero-sequence from the
voltage signals.  This is common practice, anyway, for the
zero-sequence to be generally treated separately or simply just
removed from the positive and negative ones in power system
measurements and protections.

Next, we illustrate \eqref{eq:3ph:voltage} through the following
examples.
\begin{align*}
  \rm{E0:} \quad & V_a=V_b=V_c=12~{\rm V} \, , \\
                 & \th{bo} = - \th{co}= - {2\pi}/{3} ~{\rm rad} \, . \\
  \rm{E1:} \quad & \V{a}=12~{\rm V} \, , \ \
                   \V{b}=8~{\rm V} \, , \ \ 
                   \V{c}=12~{\rm V} \, , \\
                 & \th{bo} = - \th{co}= - {2\pi}/{3} ~{\rm rad} \, . \\
  \rm{E2:} \quad & \V{a}=\V{b}=\V{c}=12~{\rm V} \, , \\
                 & \th{bo} = - {2\pi}/{3}~{\rm rad} \, , \ \
                   \th{co}= {1.5\pi}/{3}~{\rm rad} \, . 
\end{align*}

Figure \ref{fig:E0E1E2} shows the phase voltages as well as the
symmetric component of the geometric frequency ($\rho$), the Euclidean
norm of the vector component ($\omega$), and the torsional frequency
($\xi$) for examples E0 to E2.
Figure \ref{fig:vavbvc:E0E1E2} illustrates the curve formed by the
three-phase voltage in the space $(v_a, v_b, v_c)$.  Notice that
$\xi=0$ holds for the three cases.  This is consistent with the fact
that each curve in Fig.~\ref{fig:vavbvc:E0E1E2} lies in a plane.

\begin{figure*}[ht!]
  \centering
  \resizebox{\linewidth}{!}{\includegraphics{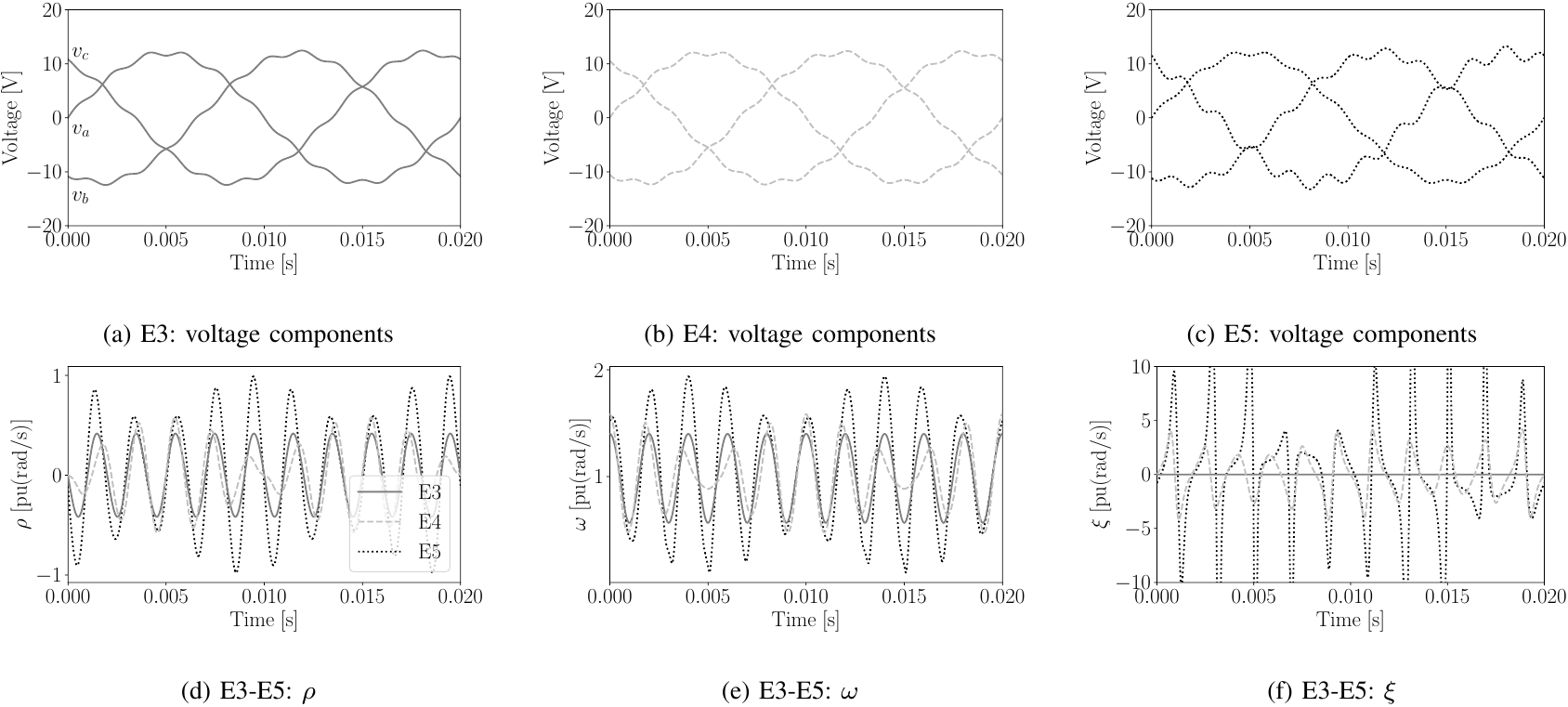}}
  \caption{Three-phase voltage components, and geometric invariants
    $\rho$, $\omega$ and $\xi$, E3-E5.}
  \label{fig:E3E4E5}
\end{figure*}

\subsection{Three-Phase AC Voltages with Harmonics}
\label{sub:harmonics}

In this section, we consider the case of three-phase voltages with
harmonics.  The conventional Fourier analysis define a basis that
consists of as many dimensions as harmonics.  While the proposed
approach can be also utilized in an arbitrary $n$-dimensional space,
we illustrate the consequences of the proposed approach in the same
three-dimensional space we have utilized in the previous section.

Let the voltage vector in \eqref{eq:3ph:vector} also include a
harmonic voltage component.  Then, the three-phase voltage becomes:
\begin{equation}
  \begin{aligned}
    \label{eq:3ph:voltharm}
    v_a &= V_a {\sin}(\wo t + \theta_{ao} )
    + V_{a,h} {\sin}(h\wo t + \theta_{ao,h})
    \, , \\ 
    v_b &= V_b {\sin}(\wo t + \theta_{bo})
    + V_{b,h} {\sin}(h\wo t + \theta_{bo,h})
    \, , \\ 
    v_c &= V_c {\sin}(\wo t + \theta_{co})
    + V_{c,h} {\sin}(h\wo t + \theta_{co,h})
    \, ,
  \end{aligned}
\end{equation} 
where $V_{i,h}$ is the magnitude of the $h$-th harmonic of phase $i$;
$\theta_{io,h}$ is the angle of the $h$-th harmonic of phase $i$ at
$t=0$~s.  Let $V_a=V_b=V_c=12$~V,
$\theta_{bo}= - \theta_{co}= - {2\pi}/{3} $~rad.  For the sake of
example, we consider the following cases of \eqref{eq:3ph:voltharm}:
E3, balanced voltage; E4, phase imbalance in the harmonic; E5,
magnitude imbalance in the harmonic.  The following values are used.
\begin{align*}
  {\rm E3} : \quad & \V{a,11}=\V{b,11}=\V{c,11}=0.5~{\rm V} \, , \\
                   & \th{bo,11} = - \th{co,11}= - {2\pi}/{3} ~{\rm rad} \, . \\
  {\rm E4} : \quad & \V{a,11}=\V{b,11}=\V{c,11}=0.5~{\rm V} \, , \\
                   & \th{bo,11} = -{2.7\pi}/{3} ~{\rm rad} \, \ \ 
                     \th{co,11}= {2.7\pi}/{3}~{\rm rad} \, . \\
  {\rm E5} : \quad
                   & \V{a,11}=0.5~{\rm V} \, , \ \ 
                     \V{b,11}=0.9~{\rm V} \, , \ \ 
                     \V{c,11}=1.3~{\rm V} \, , \\
                   & \th{bo,11} = -\th{co,11}= - {2\pi}/{3}~{\rm rad} \, .
\end{align*}

\begin{figure}[ht!]
  \centering
  \resizebox{0.7\linewidth}{!}{\includegraphics{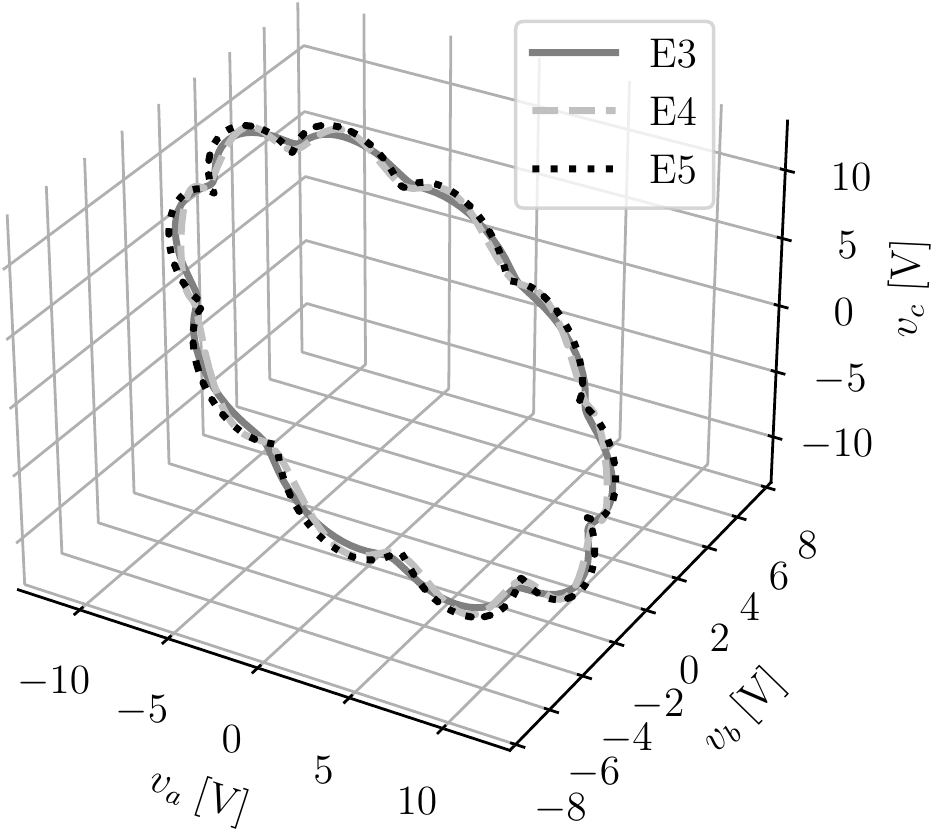}}
  \caption{Three-phase voltage in the $(v_a, v_b, v_c)$ space, E3-E5.}
  \label{fig:vavbvc:E3E4E5}
\end{figure}

The phase voltages for cases E3-E5 as well as the trajectories of
$\rho$, $\omega$ and $\xi$ are shown in Fig.~\ref{fig:E3E4E5}, while
the curve that the vector $\bfg v$ forms in the $(v_a, v_b, v_b)$
space is illustrated for each case in Fig.~\ref{fig:vavbvc:E3E4E5}.
For the balanced case, i.e.~E3, a plane curve is obtained which
implies a null torsion, whereas the curves in E4 and E5 are
three-dimensional and thus the torsion in these cases is non-zero.  It
is relevant to note that the proposed approach, differently from
Fourier analysis or any other approach based on the projections of the
signal on a kernel function (e.g., the periodic small-signal stability
analysis \cite{8113517}, wavelets and the Hilbert-Huang transform
\cite{Huang:2014}), does not need to increase the size of the basis to
take into account harmonics.

\begin{figure*}[ht!]
  \centering
  \resizebox{\linewidth}{!}{\includegraphics{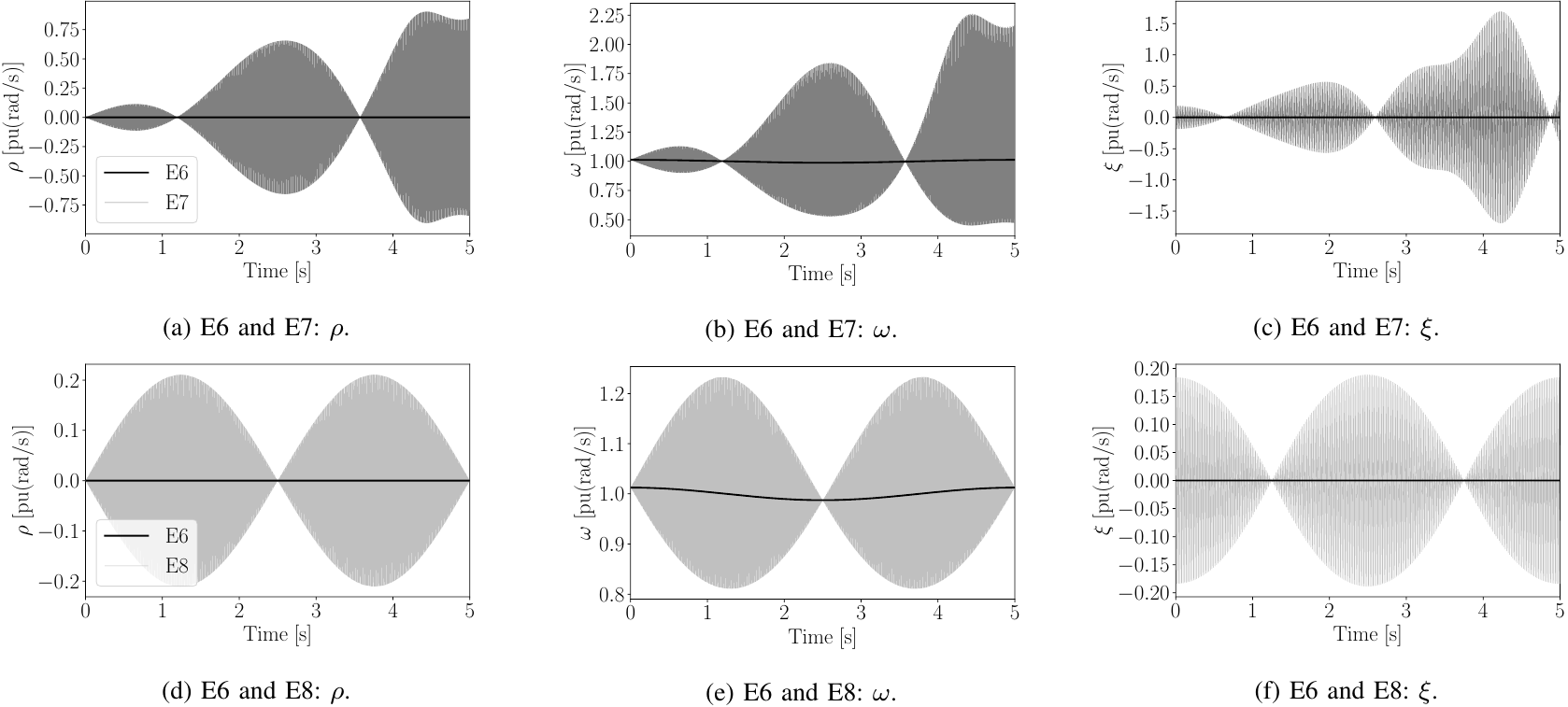}}
  \caption{Geometric frequency components and torsional frequency, E6-E8.}
  \label{fig:rhoomegaxi:E6E7E8}
\end{figure*}

\subsection{Three-Phase AC Voltages with Time-Variant Angular Frequency}
\label{sub:variable}

In the previous section, we have discussed the difference between the
Fourier approach and the proposed geometric approach.  The latter
utilizes a space that has always same dimensions regardless the number
of harmonics present in the voltage.  In this section, we further
illustrate the benefit of this frugality of dimensions.  We consider
in fact a three-phase voltage with time-varying angular frequency.  In
this case, the Fourier transform would require a basis with infinitely
many dimensions as the angular frequency varies continuously.  The
geometric approach, on the other hand, retains the three dimensions.

Let us consider that each component of the voltage vector in
\eqref{eq:3ph:vector} has a time-varying frequency $w_i$,
$i=\{a,b,c\}$.
Then, the three AC voltage components are:
\begin{equation}
  \begin{aligned}
    \label{eq:3ph:voltage:fvar}
    v_a &= V_a \, {\sin}(\w{o} t + \th{a}(t) + \th{ao}) \, , \\ 
    v_b &= V_b \, {\sin}(\w{o} t + \th{b}(t) + \th{bo}) \, , \\ 
    v_c &= V_c \, {\sin}(\w{o} t + \th{c}(t) + \th{co}) \, .
  \end{aligned}
\end{equation} 
Let $V_a=V_b=V_c=12$~V, $\w{o} = 100 \pi$ rad/s,
$\th{bo}= - \th{co}= - {2\pi}/{3} $~rad.  As an example, we consider
the following cases of \eqref{eq:3ph:voltage:fvar}: E6, balanced
$\th{i}$; E7, frequency imbalance in $\th{i}$; and E8, magnitude
imbalance in $\th{i}$.  The following values are assumed.
\begin{align*}
  {\rm E6} : \quad & \th{a} = \th{b} = \th{c} = \pi\sin(0.4\pi t) \, . \\
  {\rm E7} : \quad & \th{a} = \th{b} = \pi\sin(0.4\pi t)  \, , \\
                   & \th{c} = \pi\sin(0.44\pi t) \, . \\
  {\rm E8} : \quad & \th{a} = \th{b} = \pi \sin(0.4\pi t) \, , \\
                   & \th{c} = 1.1 \pi \sin(0.4\pi t) \, .
\end{align*}
Example~E6 is representative of the transient following a contingency
in a power system, where the oscillations of the phase angles of the
voltages are due to the electro-mechanical swings of the synchronous
machines.  On the other hand, examples E7 and E8 do not represent a
situation that can occur in a power system but are nevertheless
relevant to show the effects of the torsional frequency.

\begin{figure}[b!]
  \centering
  \resizebox{0.7\linewidth}{!}{\includegraphics{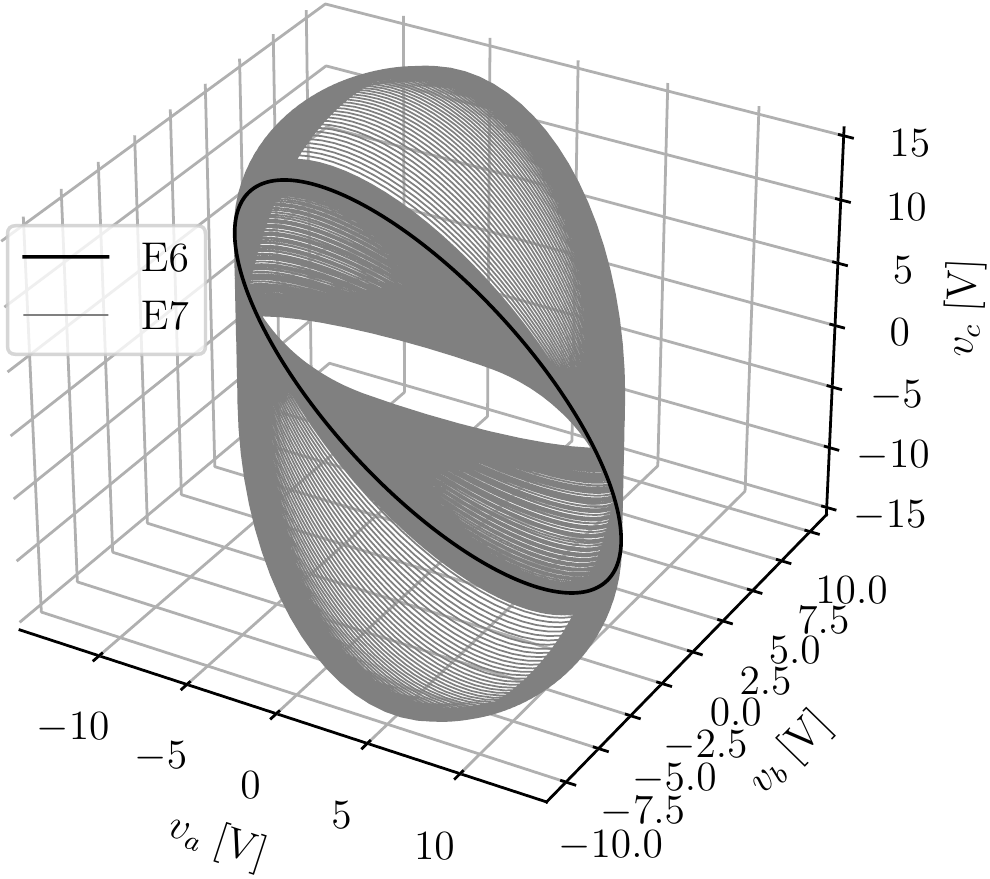}}
  \caption{Three-phase voltage in $(v_a, v_b, v_c)$ space, E6 and E7.}
  \label{fig:vabc:E7}
\end{figure}

\begin{figure}[b!]
  \centering
  \resizebox{0.7\linewidth}{!}{\includegraphics{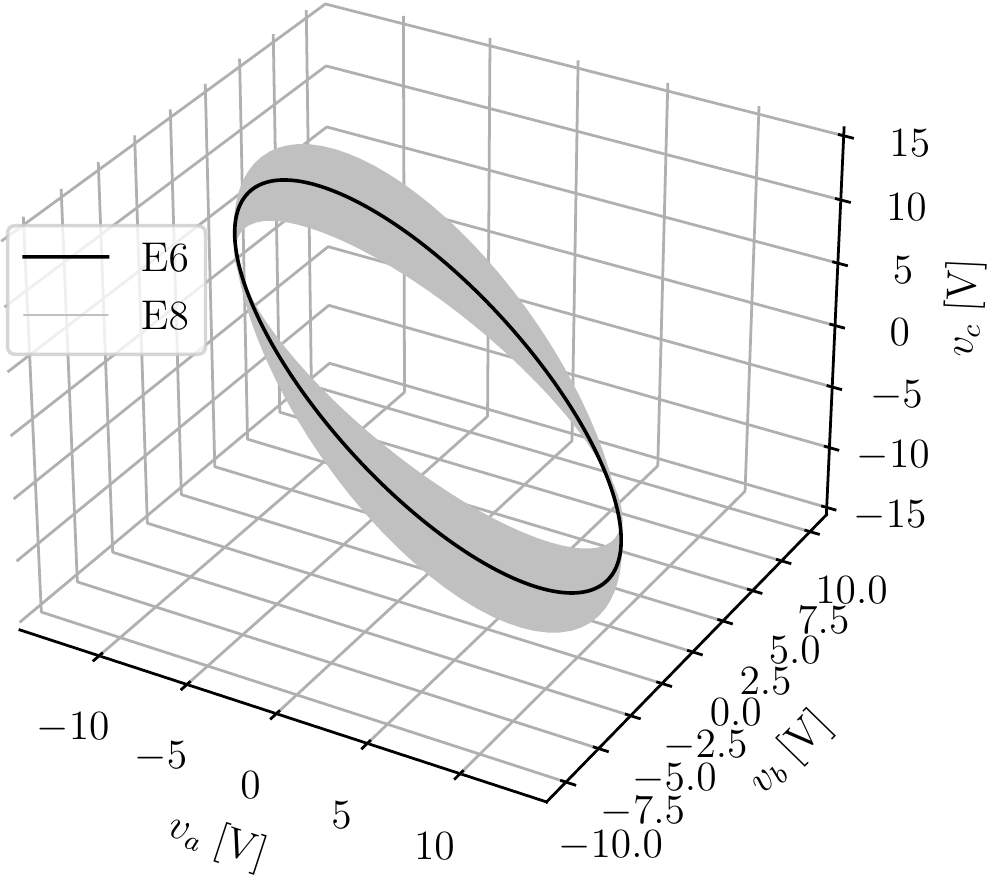}}
  \caption{Three-phase voltage in $(v_a, v_b, v_c)$ space, E6 and E8.}
  \label{fig:vabc:E8}
\end{figure}

Figure~\ref{fig:rhoomegaxi:E6E7E8} shows $\rho$, $\omega$ and $\xi$
for cases E6-E8.  In the balanced case, i.e.~E6, the three-phase
voltages show the same frequency variation and null $\rho$ and $\xi$.
The three-dimensional representation of $\bfg v$ in the space
$v_a$-$v_b$-$v_c$ is shown in Figs.~\ref{fig:vabc:E7} and
\ref{fig:vabc:E8}.  For E7 and E8, the period of the voltage is
varying with time yet not in the same way in all three phases, thus
leading to the curves shown in Figs.~\ref{fig:vabc:E7} and
\ref{fig:vabc:E8}.

\begin{figure*}[ht!]
  \centering
  \resizebox{\linewidth}{!}{\includegraphics{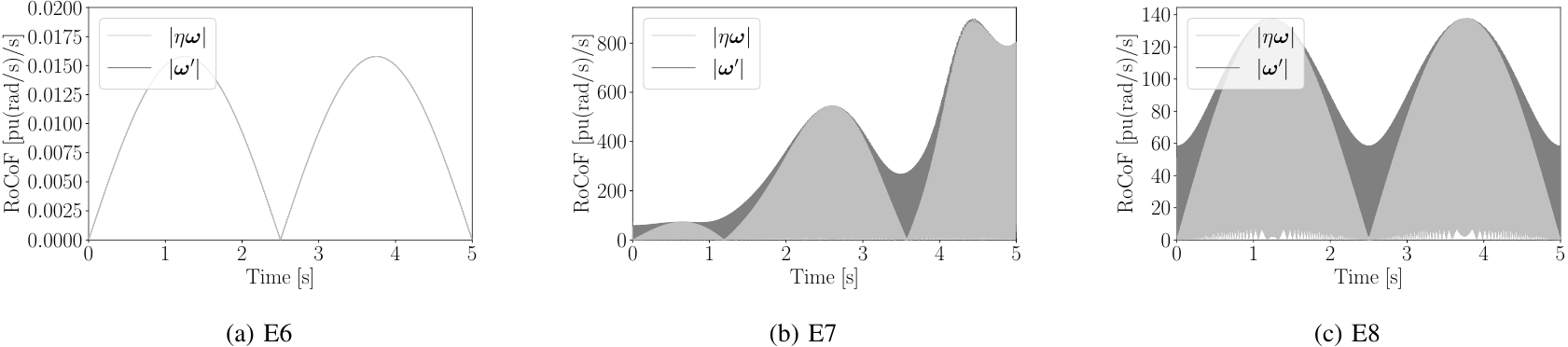}}
  \caption{Magnitudes of $\sw'$ and its symmetrical component
    $\omega' = |\eta \sw|$, E6-E8.}
  \label{fig:rocof}
\end{figure*}

Figure~\ref{fig:rocof} shows how the magnitude of the rate of change
of the vector frequency ($|\sw'|$) compares to the magnitude of its
symmetric component, i.e.~$|\eta{\sw}|$ (see equation
\eqref{eq:rocof}) for the three examples E6-E8.  These two quantities
are equal only at the time instants when the torsion is null or,
equivalently, when $\xi=0$ (see Fig.~\ref{fig:rhoomegaxi:E6E7E8}f).
It is also interesting to observe how unbalanced harmonics lead to
large variations of $|\sw'|$.


\subsection{Park Transform}
\label{sub:park}

In this final example, we consider the time derivative of the voltage
in the Park reference frame.  Let us consider a voltage vector in the
$\rm dqo$ coordinates:
\begin{equation}
  \bfg v = \vd \, \e{d} + \vq \, \e{q} + \vo \, \e{o} \, .   
\end{equation}
The time derivative of this vector in the ``inertial'' reference is
given by:
\begin{equation}
   \label{eq:vdotpark}
   \bfp v = (\vdp - \wdq \vq) \, \e{d} +
   (\vqp  + \wdq \vd) \, \e{q} + \vop \, \e{o} \, ,
\end{equation}
where $\wdq$ is the angular speed of the Park reference frame and we
have assumed a choice of the $\rm dqo$-axis such that
$\ep{d} = \wdq \e{q}$, $\ep{q} = -\wdq \e{d}$, and $\ep{o} = \bfg 0$
\cite{FDF:2020}.

The time derivative in the inertial reference in \eqref{eq:vdotpark}
is composed of two parts, namely the derivative in the rotating
$\rm dqo$-axis frame plus the effect of rotation:
\begin{equation}
  \label{eq:vdotpark2}
  \bfp v = \vrp + \bfg r \times \bfg v \, ,
\end{equation}
where:
\begin{equation}
  \vrp =  \vdp \, \e{d} + \vqp  \, \e{q} + \vop \, \e{o}  \, ,
\end{equation}
and
\begin{equation}
  \bfg r = \wdq \, \e{o} \, .
\end{equation}
Next, we apply \eqref{eq:vdot} and show that the proposed approach
leads to the same results as \eqref{eq:vdotpark} but with a different
structure of the components of $\bfp v$ with respect to
\eqref{eq:vdotpark2}.  The quantities $\rho$ and $\bfg \omega$ are:
\begin{equation}
  \rho = \frac{\vd \vdp + \vq \vqp + \vo \vop}{v^2} \, ,
\end{equation}
and
\begin{equation}
\begin{aligned}
  \sw = &\frac{\vq \vop - \vo \vqp - \wdq \vo \vd}{v^2} \, \e{d} + \\
  &\frac{\vo \vdp - \vd \vop - \wdq \vo \vq}{v^2} \, \e{q} + \\ 
  &\frac{\vd \vqp - \vq \vdp+ \wdq (\vd^2 + \vq^2)}{v^2} \, \e{o} \, ,
  \end{aligned}
\end{equation}
where $v^2 = \vd^2 + \vq^2 + v_o^2$.
Although requiring some tedious algebraic manipulations, it is not
difficult to show that, in fact, $\rho \bfg v + \sw \times \bfg v$ is
equal to the right-hand side of \eqref{eq:vdotpark}.  On the other
hand, it is straightforward to observe that:
\begin{equation}
  \vrp \ne \rho \bfg v \, , \qquad
  \bfg r \times \bfg v \ne \sw \times \bfg v \, ,
\end{equation}
i.e., $\vrp$ is not equal to the symmetric component of the time
derivative and the term $\bfg r \times \bfg v$ is not equal to the
antisymmetric component of the time derivative.

A relevant case is for $v_o = 0$, namely balanced conditions, which
leads to:
\begin{equation}
    \label{eq:park}
    \rho = \frac{\vd\vd' + \vq\vq'}{v^2} \, , \quad
    \sw = \frac{\vd \vqp - \vq \vdp + \wdq v^2}{v^2} \, \e{o} \, , 
\end{equation}
where $v^2 = \vd^2 + \vq^2$.  We note that $\rho = v'/v$, as expected,
and:
\begin{equation}
  \frac{\vd \vqp - \vdp \vq}{v^2} =  \frac{d}{dt}
  \arctan \left ( \frac{\vq}{\vd} \right ) = \Delta \omega \, ,
\end{equation}
which is the deviation of the angular frequency of $\bfg v$ with
respect to $\wdq$.  Hence the inertial- and rotating-frame time
derivatives of the voltage can be written as:
\begin{equation}
  \label{eq:park2}
  \bfp v = \rho \bfg v + (\wdq + \Delta \omega) \, \e{o} \times \bfg v \, ,
\end{equation}
and:
\begin{equation}
  \label{eq:vdotr2}
  \vrp = \rho \bfg v + \Delta \omega \, \e{o} \times \bfg v \, .
\end{equation}
The following remarks are relevant.
\paragraph*{Remark 5} The only case for which $\rho \bfg v = \vrp$ and
$\bfg r \times \bfg v = \sw \times \bfg v$ is if the angular speed of
the Park transform is the actual frequency of $\bfg v$.  In this case,
in fact, $\Delta \omega = 0$.  This result is consistent with the
commonly-used Park-Concordia model of the synchronous machine
\cite{Sauer:1998}.
\paragraph*{Remark 6} The Clarke transform can be viewed as a special
case of the Park transform, with $\wdq = 0$.  In this case
\eqref{eq:vdotpark2} leads to $\bfp v = \vrp$.  This result is
consistent with \eqref{eq:park2}.
\paragraph*{Remark 7} In stationary conditions, i.e., for
$\vdp = \vqp = 0$, \eqref{eq:park2} leads to the same expression as
\eqref{eq:vdotac}, thus confirming that balanced stationary
three-phase voltages are equivalent to a single-phase phasor.

\subsection{IEEE 39-Bus System}

This example shows the application of the proposed formulas to
estimate $\omega$, $\rho$, and $\xi$ at a bus of a power system
following a contingency.  To this aim, we use the model of the IEEE
39-bus system for Electro-Magnetic Transient (EMT) simulations
provided by DIgSILENT PowerFactory.  The system model is based on the
original IEEE 39-bus benchmark network, which has been modified to
capture the behavior during EMTs of the power network, namely, the
frequency dependency of transmission lines and the non-linear
saturation of transformers.

\begin{figure}[ht!]
  \begin{center}
    \resizebox{1.0\linewidth}{!}{\includegraphics{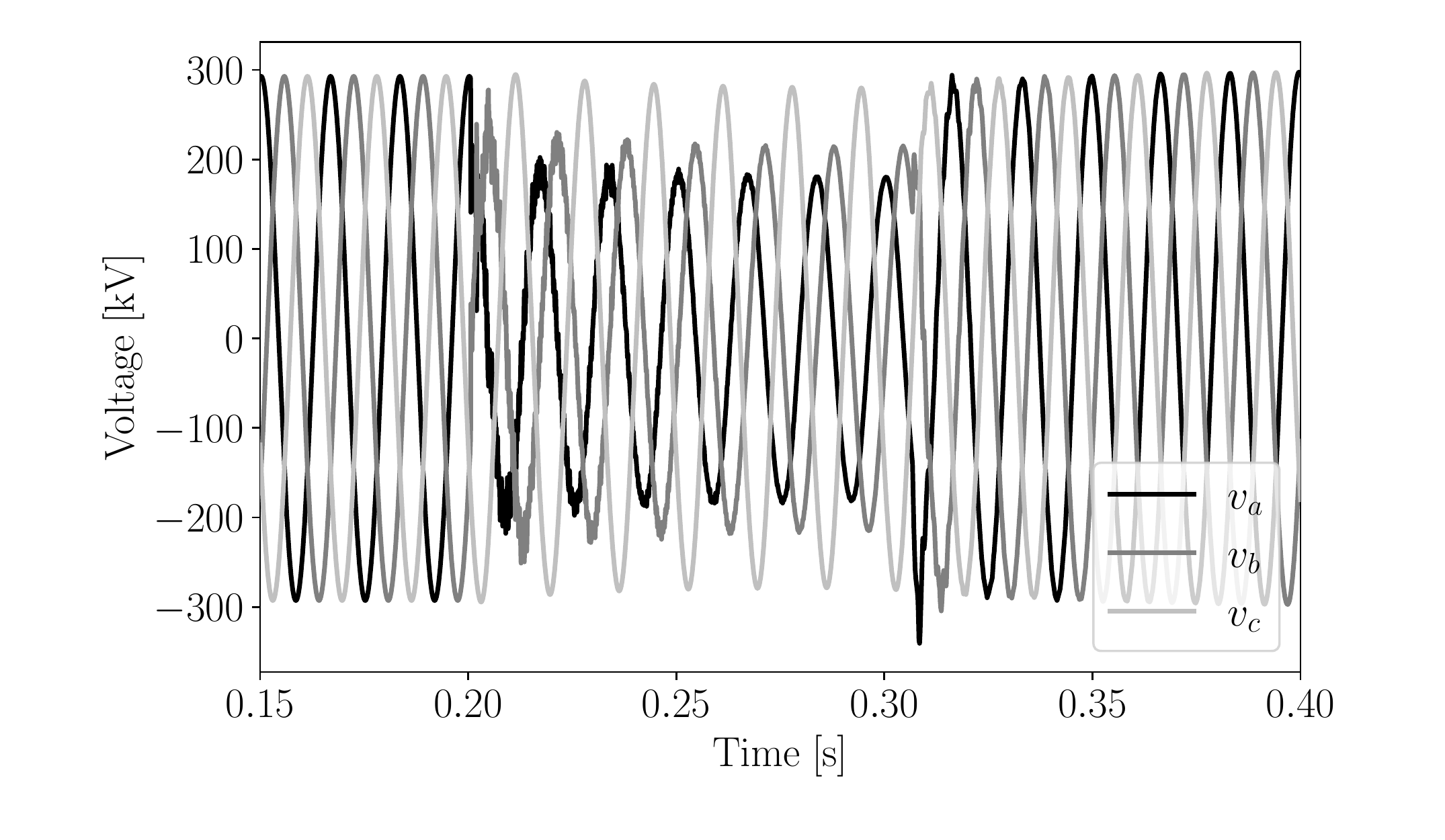}}
    \caption{Three-phase voltage at bus~26, IEEE 39-bus system.}
    \label{fig:volt}
  \end{center}
\end{figure}

\begin{figure}[ht!]
  \begin{center}
    \resizebox{0.8\linewidth}{!}{\includegraphics{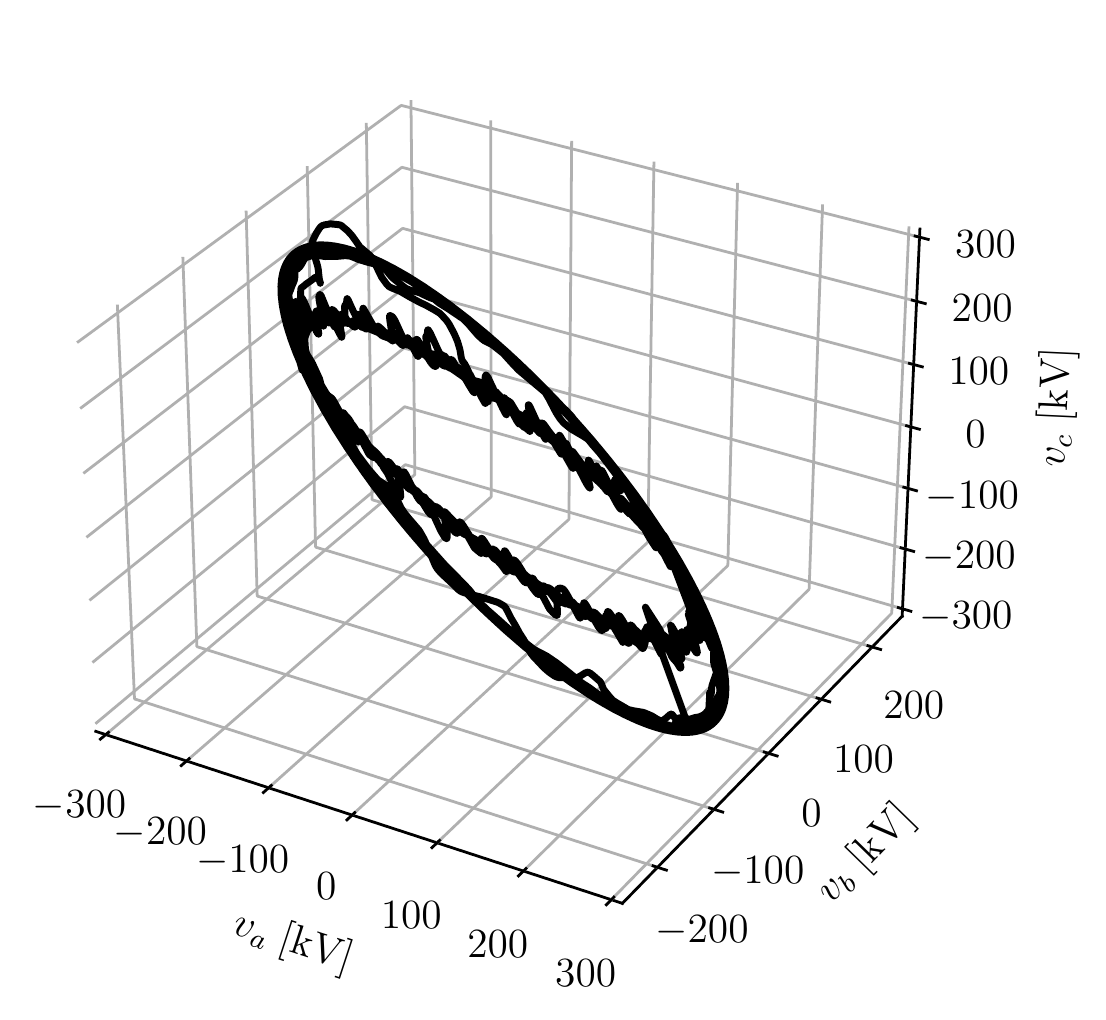}}
    \caption{Three-phase voltage at bus~26 in the space
      $(v_a, v_b, v_c)$, IEEE 39-bus system.}
    \label{fig:volt3d}
  \end{center}
\end{figure}

\begin{figure*}[ht!]
  \centering
  \resizebox{\linewidth}{!}{\includegraphics{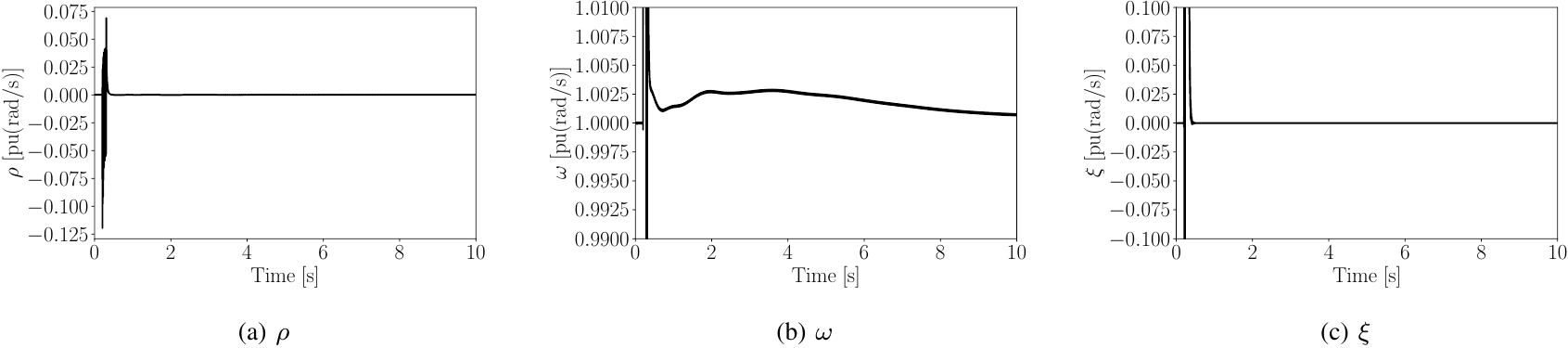}}
  \caption{Geometric invariants $\rho$, $\omega$ and $\xi$ at bus~26,
    IEEE 39-bus system.}
  \label{fig:rhowxi:39bus}
\end{figure*}

The system is numerically integrated assuming a phase-to-phase fault
between phases $a$ and $b$ at terminal bus~3 of the system at
$t = 0.2$~s. The fault is cleared at $t = 0.3$~s.  The integration
time step considered is $10^{-5}$~s.
The phase voltages at bus~26 following the contingency are shown in
Fig.~\ref{fig:volt}, while the curve formed by the three-phase voltage
in the space $(v_a, v_b, v_c)$ is illustrated in
Fig.~\ref{fig:volt3d}.  Figure~\ref{fig:rhowxi:39bus} shows $\omega$,
$\rho$, and $\xi$ following the contingency, where first-order
filtering has been applied to smooth the numerical noise in the
calculation of the voltage vector time derivatives.
Before the occurrence of the fault, the three phases are balanced and
thus, the corresponding part of the curve in Fig.~\ref{fig:volt3d} is
circular and lies in a plane.  The same holds after the fault
clearance.  Results also indicate that before the occurrence and after
the clearance of the fault, both $\rho$ and $\xi$ are null, which is
consistent to the discussion of Section~\ref{sub:variable} (e.g.,
example E6).  On the other hand, the voltage phases are unbalanced
during the fault, which gives rise to the non-circular and non-planar
sections observed in Fig.~\ref{fig:volt3d}.  For this part, both
$\rho$ and $\xi$ are non-zero as shown in Fig.~\ref{fig:rhowxi:39bus}.
This is again consistent with Section~\ref{sub:variable} and in
particular with the discussion of example E8.  Furthermore, $\omega$
accurately captures the primary frequency response at bus~26 following
the contingency (see~Fig.~\ref{fig:rhowxi:39bus}b).

Finally, for the sake of comparison, we mention that the differences
of the IEEE 39-bus system with respect to the results of E8 are that
(i) the frequency oscillation is damped and reaches a new steady state
condition, whereas this does not hold for E8, and (ii) the imbalance
occurs only for few voltage cycles, whereas in E8 voltages are
unbalanced during the whole simulation, which is the reason why the
curve in Fig.~\ref{fig:volt3d} does not appear like a compact
three-dimensional object as is the case in Fig.~\ref{fig:vabc:E8}.


\section{Conclusions}
\label{sec:conc}

The paper elaborates on the geometrical interpretation of electric
quantities and deduces several expressions that link the time
derivatives of the voltage, current and frequency in electrical
circuits with the Frenet frame.  Among these expressions, we mention
in particular \eqref{eq:vdot} and \eqref{eq:rocof2}.  Equation
\eqref{eq:vdot} indicates that the time derivative of the voltage (and
the current) is composed of two parts, one symmetric, that depends
only on the magnitude, and one antisymmetric that depends on the
``rotation'' of the quantity itself.  Equation \eqref{eq:rocof2} shows
that the time derivative of the vector frequency is more complex than
the common notion of \ac{rocof} and includes a ``rotational'' and a
``torsional'' component.  The latter is defined in this paper for the
first time.  It is interesting to note that the antisymmetric
component of the \ac{rocof} may affect the implementation and/or
performance of existing controllers.  Since the proposed approach
allows separating the symmetric and antisymmetric terms, it appears as
a useful tool for the study of power system transients and the design
of controllers.  More in general, we believe that the proposed
approach may find relevant applications in estimation, control and
stability analysis of power systems.

The proposed theory is certainly more complex than the current
conventional approach based on phasors.  However, it shows added
values from the theoretical point of view, as follows.
\begin{itemize}
\item It is a generalization of the conventional approach. The
  conventional approach, in fact, appears to be a special case of the
  proposed theory.
\item It is an example of interdisciplinary approach.  Differential
  geometry and the Frenet frame, in fact, were originally developed
  for mechanical systems.  Their applications, under certain
  hypotheses, to electrical circuits appears as an interesting advance
  which paves the way to several further developments.
\end{itemize}

An interesting byproduct of the latter point is that the proposed
theory allows ``visualizing'' electrical quantities.  This is
important, as, in the experience of the first author, students always
struggle with the lack of visual aid when studying circuit theory.
Such a support is a given in mechanical engineering.  Thus, the
ability to re-utilize well-known concepts such as curvature and torsion
also adds a didactic value to the proposed approach.

We anticipate several future work directions.  Among these, we mention
the development of a geometric framework for circuit analysis; the
applications of the formulas to estimate unbalanced conditions in
three-phase circuits as well as to circuits with more than three
phases using Cartan's extensions of the Frenet framework (see, e.g.,
\cite{Griffiths:1974}); and the development of active controllers to
reduce the effect of harmonics and imbalances.

\appendix

In this appendix, we prove the identity $n = \omega v$.  From
\eqref{eq:rhoOmega} on has:
\begin{equation}
  \label{eq:a1}
  \omega = \frac{|\sw|}{v^2} = \frac{|\bfg v \times \bfp v|}{v^2} \, .
\end{equation}
Let us focus on the term $|\bfg v \times \bfp v|$.  This can be
written as:
\begin{equation}
  \label{eq:a2}
  |\bfg v \times \bfp v| = \sqrt{(\bfg v \times \bfp v)
    \cdot (\bfg v \times \bfp v)} \, .
\end{equation}
From the following identity of the triple scalar product:
\begin{equation}
  \bfg a \cdot \bfg b \times \bfg c = \bfg b \cdot \bfg c \times \bfg a \, , 
\end{equation}
equation \eqref{eq:a2} can be rewritten as:
\begin{equation}
  \label{eq:a3}
  |\bfg v \times \bfp v| = \sqrt{\bfg v \cdot \bfp v
    \times (\bfg v \times \bfp v)} \, .
\end{equation}
Then, from the following identity of the triple vector product:
\begin{equation}
  (\bfg {a} \times \bfg {b}) \times \bfg {c} =
  (\bfg {a} \cdot \bfg {c} ) \, \bfg {b} -
  (\bfg {b} \cdot \bfg {c} ) \, \bfg {a} \, ,
\end{equation}
equation \eqref{eq:a3} can be rewritten as:
\begin{equation}
  \label{eq:a4}
  |\bfg v \times \bfp v| =
  \sqrt{\bfg v \cdot [(\bfp v \cdot \bfp v) \bfg v -
    (\bfp v \cdot \bfg v) \bfp v]} \, .
\end{equation}
From \eqref{eq:rhoOmega} and \eqref{eq:rho}, the previous expression
is equivalent to:
\begin{equation}
  \label{eq:a5}
  \begin{aligned}
    |\bfg v \times \bfp v| &= \sqrt{(\bfp v \cdot \bfp v)(\bfg v \cdot \bfg v)
      - (\bfg v \cdot \bfp v)^2} \\
    &= \sqrt{|\bfp v|^2 v^2 - \rho^2 v^4} \, ,
  \end{aligned}
\end{equation}
and, hence, \eqref{eq:a1} becomes:
\begin{equation}
  \label{eq:a6}
  \omega = \frac{\sqrt{|\bfp v|^2 - \rho^2 v^2}}{v} \, ,
\end{equation}
which, recalling the definition of $n$ given in \eqref{eq:n},
demonstrates that $n = \omega v$ and, hence,
$\bfg n = \sw \times \bfg v$.  From this relationship and the
properties of the vectors of the Frenet frame, the following
relationships follow:
\begin{equation}
  \begin{aligned}
    \bfg v = \frac{\bfg n \times \sw}{\omega^2} \, , \qquad
    \sw = \frac{\bfg v \times \bfg n}{v^2} \, .
  \end{aligned}
\end{equation}

%


\newpage

\begin{IEEEbiography}[{\includegraphics[width=1in, height=1.25in, clip,
    keepaspectratio]{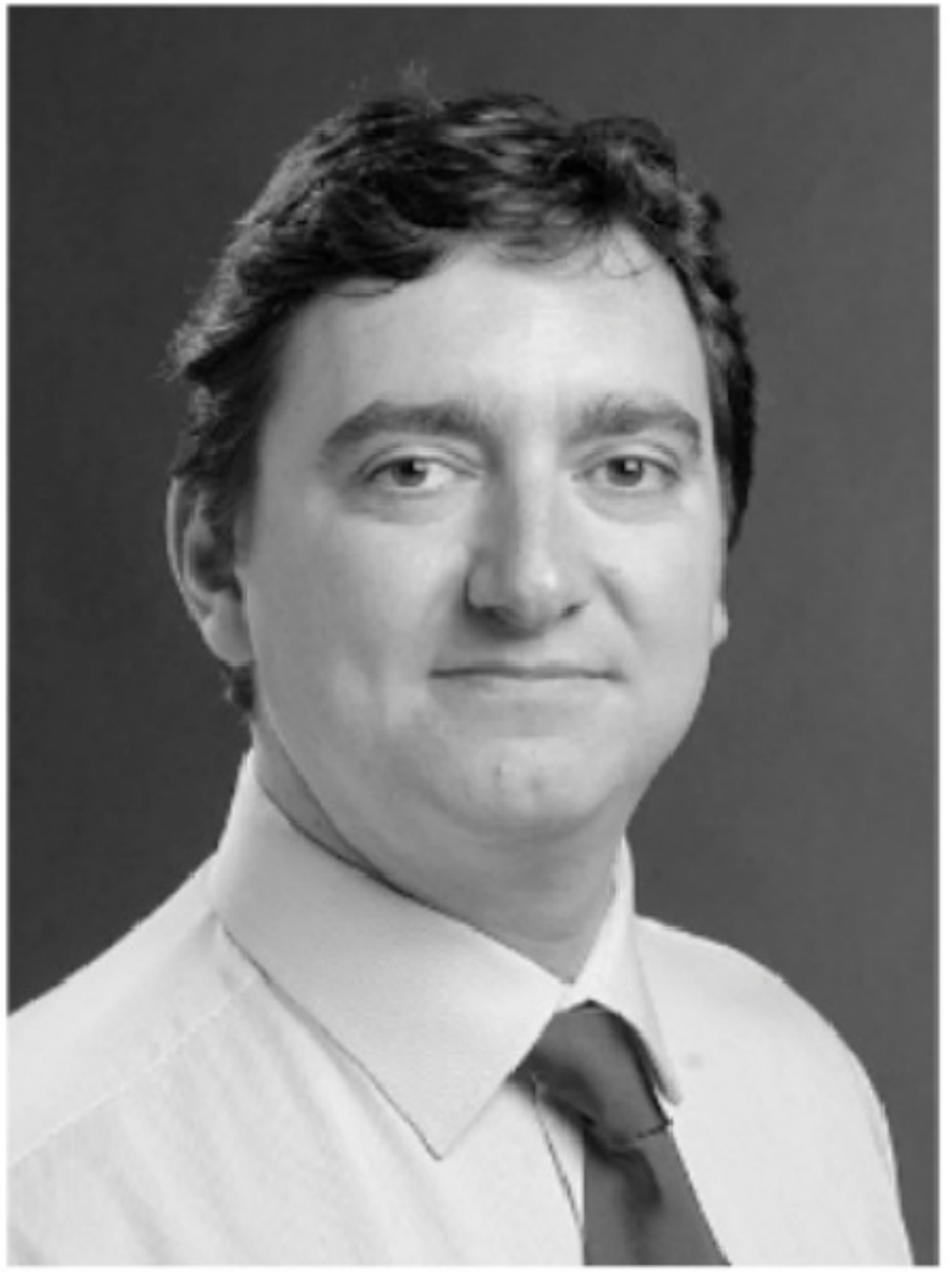}}]
  {Federico Milano} (F'16) received from the Univ. of Genoa, Italy,
  the ME and Ph.D.~in Electrical Engineering in 1999 and 2003,
  respectively.  From 2001 to 2002 he was with the University of
  Waterloo, Canada, as a Visiting Scholar.  From 2003 to 2013, he was
  with the Universtiy of Castilla-La Mancha, Spain.  In 2013, he
  joined the University College Dublin, Ireland, where he is currently
  Professor of Power Systems Protection and Control.  He is an IEEE
  PES Distinguished Lecturer, an editor of the IEEE Transactions on
  Power Systems and an IET Fellow.  He is the chair of the IEEE Power
  System Stability Controls Subcommittee.  His research interests
  include power system modelling, control and stability analysis.
\end{IEEEbiography}

\begin{IEEEbiography}[{\includegraphics[width=1in, height=1.25in, clip,
    keepaspectratio]{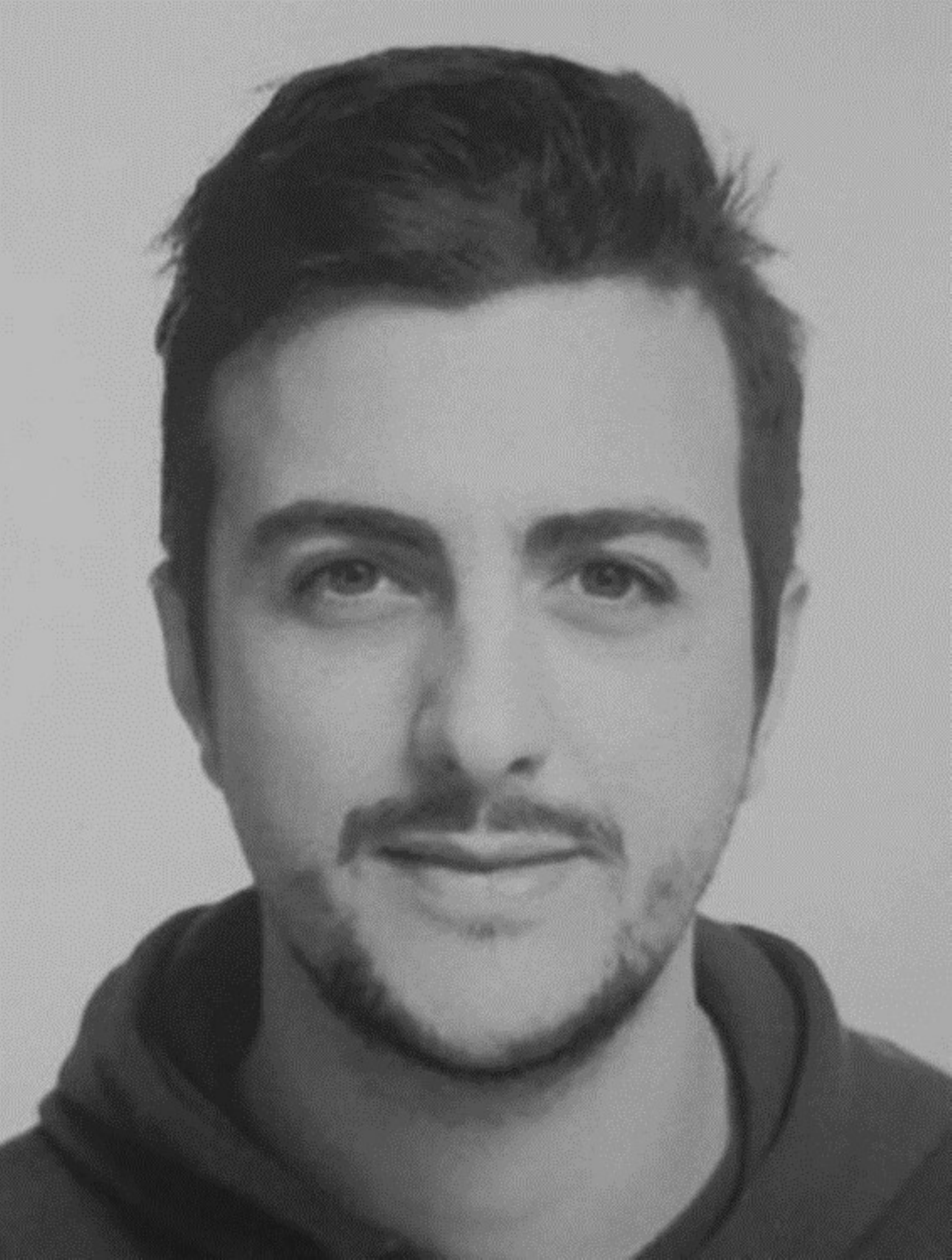}}]
  {Georgios Tzounas} (M'21) received from National Technical
  University of Athens, Greece, the Diploma (ME) in Electrical and
  Computer Engineering in 2017, and the Ph.D.~in Electrical
  Engineering from University College Dublin, Ireland, in 2021.  He is
  currently a post doctoral researcher with University~College Dublin,
  working on the Horizon 2020 project edgeFLEX. His research interests
  include modelling, stability analysis, and automatic control of
  power systems.
\end{IEEEbiography}

\begin{IEEEbiography}
  [{\includegraphics[width=1in, height=1.25in, clip,
    keepaspectratio]{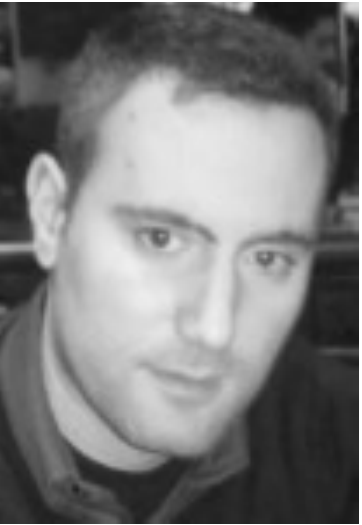}}]
  {Ioannis Dassios} received his Ph.D.  in Applied Mathematics from
  the Dpt of Mathematics, Univ.  of Athens, Greece, in 2013.  He
  worked as a Postdoctoral Research and Teaching Fellow in
  Optimization at the School of Mathematics, Univ.  of Edinburgh, UK.
  He also worked as a Research Associate at the Modelling and
  Simulation Centre, University of Manchester, UK, and as a Research
  Fellow at MACSI, Univ.  of Limerick, Ireland.  He is currently a UCD
  Research Fellow at UCD, Ireland.
\end{IEEEbiography}

\begin{IEEEbiography}
  [{\includegraphics[width=1in, height=1.25in, clip,
    keepaspectratio]{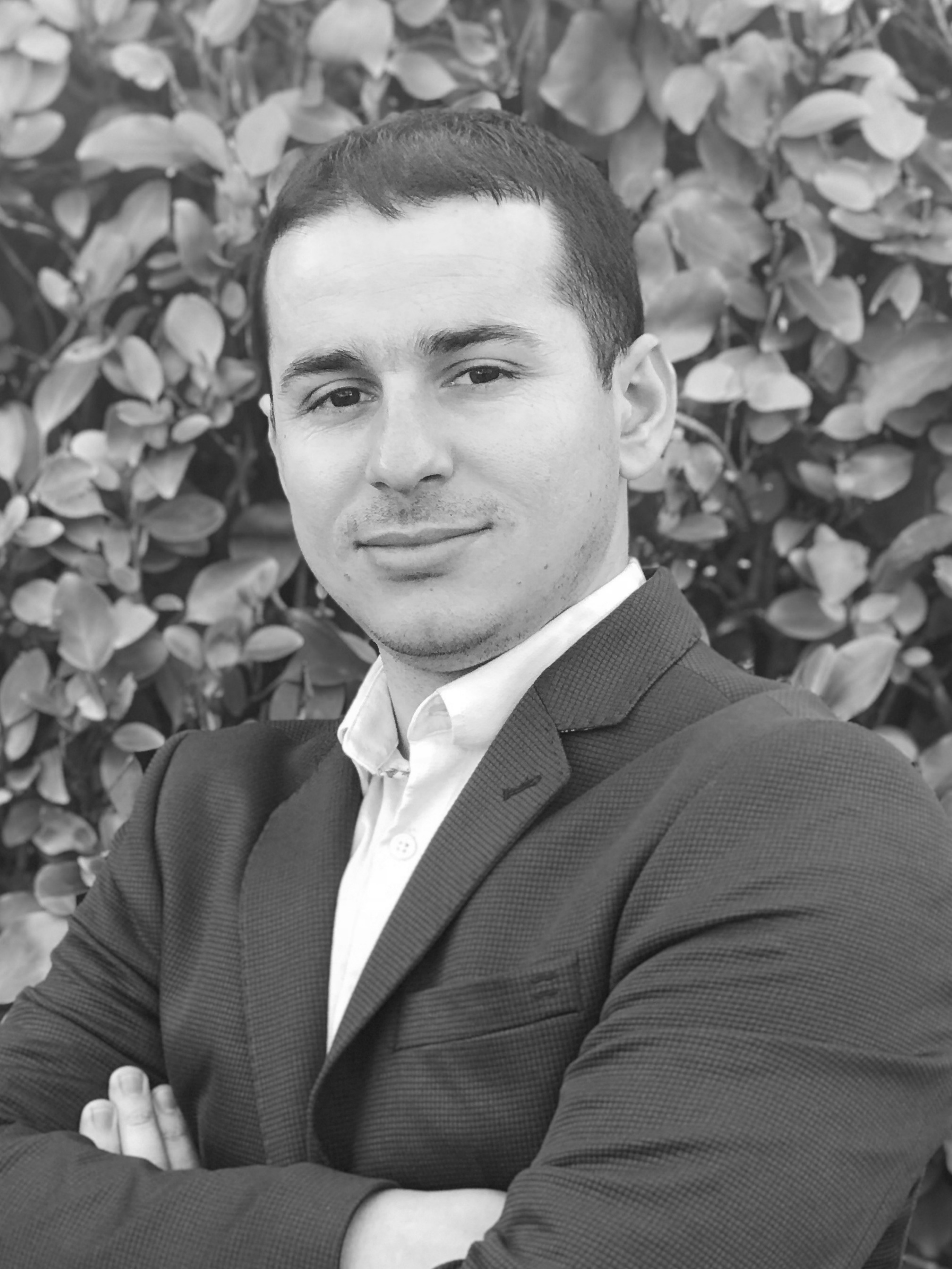}}]
  {Taulant K\"{e}r\c{c}i} (S'18) received from the Polytechnic
  University of Tirana, Albania, the BSc.~and MSc.~degree~in
  Electrical Engineering in 2011 and 2013, respectively.  From June
  2013 to October 2013, he was with the Albanian DSO at the metering
  and new connection department.  From November 2013 to January 2018,
  he was with the TSO at the SCADA/EMS office.  Since February 2018,
  he is a Ph.D.~candidate with UCD, Ireland.  In September 2021, he
  joined the Irish TSO, EirGrid.  His research interests include power
  system dynamics and co-simulation of power systems and electricity
  markets.
\end{IEEEbiography}

\vfill

\end{document}